\documentclass[12pt]{article}

\usepackage{amsmath}
\usepackage[utf8]{inputenc}

\newtheorem{theorem}{Theorem}
\newtheorem{example}{Example}
\newtheorem{lemma}{Lemma}
\newtheorem{corollary}{Corollary}

\newtheorem{proposition}{Proposition}
\newtheorem{assumption}{Assumption}
\newtheorem{algorithm}{Algorithm}
\newtheorem{definition}{Definition}

\newtheorem{rema rk}{Remark}
\newcommand{\proof}{\bf Proof: \rm }

\def\squareforqed{\hbox{\rlap{$\sqcap$}$\sqcup$}}
\def\qed{\ifmmode\else\unskip\quad\fi\squareforqed}

\def\ba{\begin{array}}
\def\ea{\end{array}}
\def\beann{\begin{eqnarray*}}
\def\eeann{\end{eqnarray*}}
\def\bea{\begin{eqnarray}}
\def\eea{\end{eqnarray}}
\def\beq{\begin{equation}}
\def\eeq{\end{equation}}
\def\Diag{\mathrm{Diag}}
\def\BT{\begin{theorem}}
\def\ET{\end{theorem}}
\def\BE{\begin{example}}
\def\EE{\end{example}}
\def\BL{\begin{lemma}}
\def\EL{\end{lemma}}
\def\BP{\begin{proposition}}
\def\EP{\end{proposition}}
\def\BC{\begin{corollary}}
\def\EC{\end{corollary}}
\def\BD{\begin{definition}}
\def\ED{\end{definition}}
\def\BA{\begin{assumption}}
\def\EA{\end{assumption}}
\def\BAL{\begin{algorithm}}
\def\EAL{\end{algorithm}}

\def\R{{\rm I\kern-3.2pt R}}
\def\by{\bar{y}}

\def\det{\mbox{det }}
\def\vol{\mathrm{vol}}

\begin{document}

\title{The ellipsoid method redux}
\author{
Michael J.~Todd
\thanks{
School of Operations Research and Information Engineering,
Cornell University, Ithaca, NY 14853,
USA.
E-mail {\tt mjt7@cornell.edu}.
}}
\maketitle

{ \small \hspace{.5 in}  Dedicated to the memory of Naum Z.\ Shor and to the Ukrainian people.}

\abstract{
We reconsider the ellipsoid method for linear inequalities. Using the
ellipsoid representation of Burrell and Todd, we show the method can be viewed as coordinate descent on
the volume of an enclosing ellipsoid, or on a potential function, or on both.
 The method can be enhanced by improving the lower bounds generated and by allowing the weights on inequalities to be
decreased as well as increased, while still guaranteeing a decrease in volume.
Three different initialization schemes are described, and preliminary computational results given. Despite the improvements 
discussed, these are not encouraging.
}



\section{Introduction}\label{sec:intro}

In the mid 1970s, a Russian, Arkadi S.\ Nemirovski, and a Ukrainian, Naum Z.\ Shor, independently
devised the ellipsoid method for convex nonsmooth minimization. Their motivations were rather different:
Yudin and Nemirovski \cite{yn} were interested in the informational complexity of convex optimization and in
developing an implementable version of the method of central sections of Levin \cite{lev} and Newman \cite{new},
while Shor \cite{shor} was investigating a special case of his space dilation methods with intriguing geometric
properties. The history of the study of convex minimization in the Soviet Union is nicely described in the survey
article of  Tikhomirov \cite{tikh}.

The algorithm did not attract a great deal of interest in the West until a couple of years later, when Leonid G.\
Khachiyan \cite{kha} used it to prove the polynomial-time solvability of linear programming. This was a very impressive
theoretical advance, but the algorithm did not seem to be useful in solving large-scale sparse linear
programming problems in practice. There were also important consequences in combinatorial optimization, as
noted by Karp and Papadimitriou \cite{kp},  Gr\"{o}tschel, Lov\`{a}sz, and Schrijver
\cite{glspaper}, and Padberg and Rao \cite{pr}, and explored in depth in \cite{gls}. 

Considerable efforts were made over the years after the ellipsoid method was first developed to improve its performance.
The fundamental problem is that the volume of the enclosing ellipsoid is cut only by a factor of about $1 - (2n)^{-1}$ at each iteration, 
where $n$ is the dimension of the problem, whereas the
method of central sections cuts the volume of an enclosing set by a constant factor. Thus variants using so-called deep cuts and parallel
cuts were developed first by Shor and Gershovich \cite{Shorg} and then rediscovered by many others, but unfortunately did not
improve practical performance by much. Nemirovski and Yudin devised algorithms  adapting to the effective dimension of
the problem \cite{nemyud77}, and recently Rodomanov and Nesterov \cite{rn} have developed a hybrid of the ellipsoid and subgradient methods
avoiding the difficulties of the former in high dimensions. Tarasov, Khachiyan, and Erlikh  \cite{tke88}
developed an inscribed ellipsoid method that decreased the volume by a constant ratio,
but at the expense of a much more complicated computation at each iteration.  A summary of the early developments can be found in \cite{bgt}.

In the domain for which it was intended, the ellipsoid method did prove quite effective for medium-sized highly nonlinear
problems (Ecker and Kupferschmid \cite{ek}), although testing by Shor indicated that other versions of his space-dilation methods
were preferable; see also the more recent theoretical analysis in Burke, Lewis, and Overton \cite{blo}.

Our interest here is in linear programming, and more specifically in systems of linear inequalities. 
As we have noted, the ellipsoid method can be very slow
when attacking problems of high dimension, even when using deep cuts. However, there is a hope of improved
performance when the algorithm is viewed a different way. Burrell and Todd \cite{bt85} showed that an enclosing ellipsoid could be
derived by combining rank-one convex quadratic inequalities obtained from the original inequalities, provided these
were two-sided. Thus the Burrell-Todd version of the ellipsoid method required lower bounds on the linear functions involved in the inequalities (we assume these
are all given by upper bounds), and these bounds were obtained from duality arguments. At each iteration, a single violated inequality
was chosen, and then the multiplier $d_j$ on the associated quadratic inequality adjusted, possibly after updating the corresponding
lower bound. This method was shown to be closely related to the deep/parallel cut ellipsoid method. One advantage of this
viewpoint is that it is sometimes possible to infer the infeasibility of the original system of linear inequalities, and obtain a
certificate of this infeasibility, whereas the original method could only do this in the case of rational data, after making perturbations, taking
an astronomical number of iterations without obtaining a feasible point, and then concluding infeasibility without producing a certificate.

Recently, Lamperski, Freund, and Todd \cite{LFT} developed the {\em oblivious ellipsoid algorithm} (OEA), which by modifying the
method above slightly is guaranteed to either find a feasible solution or prove infeasibility in a number of iterations that is
polynomial in the dimensions of the problem and the logarithm of a certain condition number of the system.
In the case of infeasibility, the proof relies on a fixed decrease in a certain potential function. However, the algorithm
operates obliviously (without knowing whether the problem is feasible or infeasible), and so it must choose its parameters carefully
to guarantee a simultaneous decrease in both the volume of the enclosing ellipsoid and the potential function. As a result,
the decrease in the volume is even slower, by a factor of about $1 - (2m)^{-1}$, where $m > n$ is the number of inequalities.
Another consequence was some counter-intuitive steps in the algorithm: the ellipsoid was not always defined by the best
lower bounds obtained on the constraint functions, and sometimes these lower bounds were actually decreased.

In this paper we propose improvements to both the OEA and the standard (deep-cut) ellipsoid 
algorithm (SEA), implemented as in Burrell-Todd \cite{bt85}. Among all
the lower bounds on the $j$th constraint function $a_j^T x$ that could be derived from ellipsoids differing from the current one
by only adjusting the weight $d_j$ and its defining lower bound $l_j$, 
we show how the ``best'' one can be obtained at negligible cost. We demonstrate how the OEA
can be adapted with no loss in theoretical guarantees to always use the best lower bounds generated. The SEA can also use this
improved bound. We also show how the SEA can be modified to use decreasing or drop steps which decrease the weight 
$d_j$ on a well-satisfied constraint, again without losing theoretical guarantees. In addition, we give a modification of the potential
function which allows some weights $d_j$ to be zero, and still allows convergence to be established in the infeasible case. We provide a
geometric interpretation of these potential functions.

Finally, we give the results of some preliminary computational tests, which show that, while the performance of the algorithms is improved
by these changes, it remains unfortunately slow, with a number of iterations growing roughly with the $1.7$-th power of the dimension.

The paper is organized as follows. In Section \ref{sec:elliprep}, we give the ellipsoid representation we employ and show
some of its properties. We discuss the improved lower bound in Section \ref{sec:bestlb}. Then, in Section \ref{sec:potfn},
we define the modified potential function and show that it can also be decreased by a constant at every iteration.
We consider steps that decrease weights $d_j$ in Section \ref{sec:decsteps}, and put all these ideas together in 
improvements of the OEA and the SEA in Section \ref{sec:improving}. The results of computational testing are given in Section \ref{sec:comp}.

\section{The ellipsoid representation and certificates of infeasibility}\label{sec:elliprep}

In this section we describe the technique Burrell and Todd \cite{bt85} used to
generate an ellipsoid containing a polyhedron and how it can be used also to certify infeasibility.
We will generally use the notation of \cite{bt85}, except that, as in \cite{LFT}, we use $\gamma_j$
to denote the semi-width of the ellipsoid in the direction of $a_j$ rather than its square, to avoid
square roots.

We would like either to find a point in
\[
P := \{ y \in \R^n: A^T y \le u \}
\]
or to prove it empty.
We assume $A$ is $n \times m$, and has rank $n$ (otherwise we can project $y$ and the columns of $A$ to a lower-dimensional space).
The $j$th column of $A$ is denoted $a_j$; without loss of generality, we suppose all columns are nonzero.

Our first step is to provide lower bounds also on the linear functions in $A^T y$. Let us assume we have an
$m \times m$ matrix $\Lambda$ and an $m$-vector $l$ with
\beq \label{eq:Lambda}
A \Lambda = - A, \quad \Lambda \geq 0, \quad l = - \Lambda^T u,
\eeq
so that $y \in P$ implies
\[
A^T y = - \Lambda^T A^T y \geq - \Lambda^T u = l,
\]
implying that $l$ is a vector of lower bounds as desired. 
Burrell and Todd \cite{bt85} explain how to obtain such a $\Lambda$ if we are given upper and lower bounds
on the components of a feasible $y$. Then
\beq \label{eq:defP}
P = \{ y \in \R^n: l \leq A^T y \leq u \}.
\eeq
We will also assume $l < u$; otherwise $P$ is empty or we can again project it to a lower-dimensional space.

We now proceed from this bounded polyhedron defined by many linear constraints (with potentially exponential complexity in vertices and faces)
to an approximating set defined by a single quadratic inequality, as follows: First we multiply
linear inequality by linear inequality, to get a (rank-one) quadratic inequality; then we aggregate these inequalities
using nonnegative weights. Thus the inequalities indexed by $j$ imply
\[
(a_j^T y - l_j) (a_j^T y - u_j) \leq 0;
\]
then, multiplying this by $d_j \geq 0$ and summing, we find
\[
y \in P \implies \sum_j d_j (a_j^T y - l_j) (a_j^T y - u_j) \leq 0.
\]
If we let $D := \Diag ( d)$, this quadratic inequality can be written
\beq\label{eq:quadineq0}
(A^T y - l)^T D (A^T y - u) \leq 0.
\eeq

We always assume that $d$ is chosen so that $ADA^T$ is positive definite. Then we can complete the
square as follows. Define
\[
r := \frac{u + l}{2}, \quad v := \frac{u - l}{2},
\]
so that $l = r - v, u = r + v$. Then set
\[
\bar y := (ADA^T)^{-1} ADr.
\]
Now simple algebraic manipulations show that the quadratic inequality (\ref{eq:quadineq0}) can be written as
\beq \label{eq:quadineq}
\ba{rcl}
(y - \by)^T (ADA^T) (y - \by) & \leq & \by^T (ADA^T) \by - l^T D u  \\
    & = & r^T D A^T (ADA^T)^{-1} A D r - l^T D u  \\
     &  = &  r^T D A^T (ADA^T)^{-1} A D r - r^T D r + v^T D v.
\ea
\eeq
If the right-hand side above is positive, this defines an ellipsoid centered at $\bar y$; if zero, 
the ellipsoid degenerates to a single point; and if negative, the inequality is
infeasible.

Let us suppose for now that the right-hand side is positive. Clearly the ellipsoid depends on the vector $d$, but since the algorithm
sometimes updates the lower bounds $l$ (while also updating $\Lambda$ to preserve the derivation (\ref{eq:Lambda})),
we denote it by either $E(d)$ or $E(d,l)$. We have shown above that
\beq \label{eq:cont}
\ba{rcl}
P & \subseteq & E(d) := E(d,l)  \\
   & :=  & \{ y \in \R^n: (y - \by)^T (ADA^T) (y - \by) \leq  r^T D A^T (ADA^T)^{-1} A D r - r^T D r + v^T D v \}.
   \ea
\eeq

Since we are temporarily assuming that the right-hand side 
\beq \label{eq:deff}
f(d) := f(d,l) := r^T D A^T (ADA^T)^{-1} A D r - r^T D r + v^T D v 
\eeq
is positive, we can scale $d$ so that it becomes 1. Then, if we write
\[
B := (ADA^T)^{-1}
\]
(we will use this notation whether $f(d,l)$ is 1 or not),
the ellipsoid can be written $\{ y \in \R^n: (y - \by)^T B^{-1} (y - \by) \leq 1 \}$, which is the traditional way to
represent ellipsoids in the ellipsoid method. Let us note now two important advantages of the representation in (\ref{eq:cont}) we
are using:
\begin{itemize}
\item[(i)] The conditions in (\ref{eq:Lambda}) and the nonnegativity of $d$ certify the containment
(\ref{eq:cont}); and

\item[(ii)] The matrix $ADA^T$ preserves more of the sparsity and structure of the matrix $A$ than its inverse $B$,
and a Cholesky factor of $ADA^T$ is also likely to preserve some of this sparsity (as in interior-point methods),
while allowing us to perform cheaply all the operations in the algorithm. In contrast, $B$ may be dense or close to dense.

\end{itemize}

Given such an ellipsoid $E$ containing $P$, an iteration of the ellipsoid method
\begin{itemize}
\item[(a)] stops if its center $\by$ satisfies all the constraints, and otherwise

\item[(b)] possibly stops with a proof of infeasibility, or

\item[(c)] generates a new ellipsoid $E_{+}$ containing $P$ and satisfying
\beq \label{eq:volred}
\vol(E_{+}) \leq \exp\left[ - \frac{1}{2(n+1)}\right] \vol(E).
\eeq

\end{itemize}

The oblivious ellipsoid algorithm (OEA) of Lamperski, Freund, and Todd \cite{LFT} is designed to
terminate in a polynomial number of steps whether the system of inequalities is feasible or not
(and without knowing which), and so takes a more conservative approach to updating the ellipsoid,
so that (\ref{eq:volred}) is modified to replace the dimension $n$ with the number of inequalities $m$:

\beq \label{eq:volred2}
\vol(E_{+}) \leq \exp\left[ - \frac{1}{2(m+1)}\right] \vol(E).
\eeq

At the same time, the OEA guarantees a decrease in a certain potential function $\phi(d,l)$ which we will define later,
so that its value at the new iterate is decreased by the same factor:

\beq \label{eq:potred}
\phi_{+} \leq \exp\left[ - \frac{1}{2(m+1)}\right] \phi.
\eeq

If the ellipsoid is represented as above, with $E= E(d,l)$ and $d$ scaled so that $f(d,l) = 1$, (b) and (c) are accomplished by first choosing a violated
constraint $j$; then checking whether $l_j$ is at least as large as the minimum value of $a_j^T y$ over $E$,
$a_j^T \by - (a_j^T B  a_j)^{1/2}$, and if not, updating it and the $j$th column of $\Lambda$ so that
(\ref{eq:Lambda}) remains true; and finally increasing just the $j$th component of $d$ so that
(\ref{eq:volred}) (or (\ref{eq:volred2}) and (\ref{eq:potred})) holds. (The OEA also modifies $l_j$ in a more complicated way; the first part of the step is to
decrease $l_j$ until the center $\by$ satisfies $a_j^T \by = u_j$.)
We will elaborate on this in the following sections; details can be found in Burrell and Todd \cite{bt85} and Lamperski et al.\ \cite{LFT}. 
We note that the algorithm in \cite{LFT} and its analysis are quite technical. We will therefore not provide more details here, hoping
that our version in Section \ref{sec:improving} is easier to understand, but we will rely on the analysis in \cite{LFT} to establish convergence.

Among the contributions of this paper are a better way to update the lower bound
and  to show that a suitable volume reduction can also be achieved by choosing a
constraint that is {\em well satisfied}, and then {\em decreasing} the corresponding component of $d$, possibly to zero. Since decreasing
a component of $d$ can make the right-hand side of the quadratic inequality negative (or zero), we are led to
another advantage of the representation we are using:
\begin{itemize}

\item[(iii)]  If $f(d,l)$ is negative (or if it is zero and the center $\by$ is not feasible),
the representation in (\ref{eq:cont}) provides a certificate of infeasibility for $P$.

\end{itemize}

Throughout this paper, a certificate of infeasibility means a Farkas-type certificate (an L-certificate in the notation of
\cite{LFT}). We say that $x \in \R^m$ certifies the infeasibility of $l \leq A^T y \leq u$ if
\beq \label{eq:certinf1}
A x = 0, \qquad u^T x_+ - l^T x_- < 0
\eeq
(here $x_+$ and $x_-$ are the componentwise positive and negative parts of $x$ so that $x_+ - x_- = x$,
$x_+ \geq 0$, $x_- \geq 0$, $\mbox{ and } x_+^T x_- = 0$). To show that such a certificate is valid, note that
 if $x$ satisfies (\ref{eq:certinf1}), and if $l \leq A^T y \leq u$, then $u^T x_+ \geq y^T A x_+$
and $l^T x_- \leq y^T A x_-$, so subtracting gives $u^T x_+ - l^T x_- \geq y^T A (x_+ - x_-) = 0$, which is impossible. 
As an example, suppose a lower bound $l_j$ is certified by a vector $\lambda_j$ as in (\ref{eq:Lambda}) so that
$A \lambda_j = - a_j, \, \lambda_j \geq 0, \, \mbox{ and }  - u^T \lambda_j = l_j$, and $l_j > u_j$; then we find
$A (\lambda_j + e_j) = 0, \, \lambda_j + e_j \geq 0, \, \mbox {and } u^T (\lambda_j + e_j) = u_j - l_j < 0$,
so that $x := \lambda_j + e_j$ is a certificate of infeasibility. Here $e_j$ denotes the $j$th coordinate vector; we also use $e$ to 
denote the vector of ones of appropriate dimension.

We note that a certificate of infeasibility for $l \leq A^Ty \leq u$ easily yields one for $A^Ty \leq u$, i.e.,
a vector $\hat x$ with
\[
A \hat x = 0, \quad \hat x \geq 0, \quad u^T \hat x < 0.
\]
Indeed, if $x$ satisfies (\ref{eq:certinf1}), then it is easy to see using (\ref{eq:Lambda}) that
$\hat x = x_+  +  \Lambda x_-$ satisfies the system above.

Item (iii) above was proved in \cite{LFT}; here we provide a simpler derivation.
We first establish the claim in the case of a negative right-hand side
(in Section \ref{sec:potfn} we will give geometric intuition for why $f(d,l) < 0$ implies
infeasibility).

\BT \label{th:negrhs}
Suppose, with the notation above, 
\[
f(d,l) = r^T D A^T B ADr - r^T D r + v^T D v < 0.\]
Then
\[
x := D \bar t ,
\]
where
\beq\label{eq:deftbar}
\bar t := A^T \by - r,
\eeq
certifies the infeasibility of $l \leq A^T y \leq u$.
\ET

\proof
We need to show that $D \bar t$ satisfies the conditions in (\ref{eq:certinf1}). . First, $Ax = 0$ follows from the definition of $\by$.
Now note that the projection of $D^{1/2} r$ onto the null space of $A D^{1/2}$ is
\beq\label{eq:defq}
q := D^{1/2} r - D^{1/2} A^T (ADA^T)^{-1} A D^{1/2} D^{1/2} r = - D^{1/2} \bar t,
\eeq
and its norm squared is
\beq\label{eq:qtq}
q^Tq = r^T D r - r^T D A^T (ADA^T)^{-1} AD r.
\eeq
So the right-hand side of the quadratic inequality being negative is equivalent to $ - q^Tq + v^T D v < 0$,
or
\[
\| D^{1/2} v \| < \| q \|.
\]
Finally, 
\[
u^T x_+ - l^T x_- = (r+v)^T x_+ - (r-v)^T x_- = r^T(x_+ - x_-) + v^T(x_+ + x_-) = r^T x + v^T |x|,
\]
where $|x|$ denotes the vector of absolute values of the components of $x$. Now $r^T x = (D^{1/2}r)^T (D^{1/2} \bar t) =
- (D^{1/2}r)^T q = - q^T q = - \| q \|^2$, while $v^T |x| = (D^{1/2} v)^T |D^{1/2} \bar t| \leq \| D^{1/2} v \| \| q \| < \| q \|^2$,
which shows that the quantity displayed above is negative.

 \indent \squareforqed

The proof incidentally provides another useful form for $f(d,l)$. Using
(\ref{eq:deff}), (\ref{eq:defq}), and (\ref{eq:qtq}), we obtain
\beq\label{eq:deff2}
f(d,l) = v^T D v - \bar t ^T D \bar t.
\eeq

Finally, we consider the case that the right-hand side $f(d,l)$ is zero. Then the solution set to the quadratic inequality
is the singleton $\by$. If this is feasible, we have our desired point in $P$. If not, we know $P$ is empty, but again
we would like a certificate of infeasibility. We show that either $x$ or a perturbation of it provides such a certificate.

\BP\label{pr:feq0}
Suppose now 
\[
f(d,l) = 0.
\]
Then, if $\by$ is infeasible, either
\[
x := D \bar t
\]
or 
\[
\tilde x := x + \epsilon (e_j - DA^T B a_j)
\]
is a certificate of infeasibility for $l \leq A^T y \leq u$, where in the second case
either $a_j^T \by > u_j$ and $\epsilon > 0$ or $a_j^T \by < l_j$ and $\epsilon < 0$,
and $\epsilon$ is sufficiently small in absolute value that none of the nonzero components of $x$ change sign in $\tilde x$.
\EP

\proof
By examining the proof of the theorem above, we see that $\|D^{1/2 } \bar t \| = \| D^{1/2} v \|$ and that, if $x$ does not
yield a certificate, $D^{1/2} | \bar t |$ and $D^{1/2} v$ are collinear. Henceforth, assume $\by$ is infeasible and $x$ is
not a certificate of infeasibility. Then, for nonzero $d_h$,  $\bar t_h  =  \pm v_h$, and 
$a_h^T \by =: w_h$ is either $l_h$ or $u_h$. Moreover, if $d_h \neq 0$ but $x_h = 0$, then $a_h^T \by = r_h$ and $\bar t_h = 0$,
so $v_h = 0$ and $w_h = l_h = u_h$. 

Thus, if $d_h$ is nonzero, $a_h^T \by = w_h$, so that $DA^T \by = Dw$. Then $ADA^T \by = ADw$, so $B ADw = \by$.
Now choose $j$ so that $a_j^T \by > u_j$ or $a_j^T \by < l_j$ (since $\by$ is
infeasible), so that $d_j = 0$, and choose $\epsilon$ as in the statement of the proposition.

Since $Ax = 0$ and $A(e_j - DA^T B a_j) = a_j - a_j$, we have $A \tilde x = 0$, and it remains to show that
$u^T \tilde x_+ - l^T \tilde x_- < 0$. If $x_h > 0$, then $w_h > r_h$, so $w_h = u_h$ and $\tilde x_h > 0, \, u_h \tilde x_h = w_h \tilde x_h$.
Similarly, if $x_h < 0$, then $\tilde x_h < 0, \,  (- l_h) (- \tilde x_h) = w_h \tilde x_h$. And if $d_h \neq 0$ but $x_h = 0$, then $w_h = u_h = l_h$,
so $u_h \tilde x_h = (- l_h)(- \tilde x_h) = w_h \tilde x_h$.  A similar argument with $x$ instead of $\tilde x$ shows
that $w^T x = u^T x_+ - l^T x_- = 0$, since $x$ is not a certificate of infeasibility. Now suppose $a_j^T \by > u_j$ and so $\epsilon > 0$. Then
\beann
u^T \tilde x_+ - l^T \tilde x_- & = & \sum_{d_h \neq 0} [w_h x_h - \epsilon (w_h d_h a_h^T B a_j )] + \epsilon u_j \\
  & = & w^T x + \epsilon( u_j - w^T D A^T B a_j)  \\
  & = & \epsilon (u_j - a_j^T \by) < 0,
 \eeann
 and a similar argument holds if $a_j^T \by < l_j$.

 \indent \squareforqed

It appears therefore that decreasing the right-hand side $f(d,l)$ aids both in finding a feasible point in $P$ if one exists
(by decreasing the volume of the enclosing ellipsoid $E(d,l)$) and also in finding a certificate of infeasibility
in the case that $P$ is empty. We now make this more precise. 

By considering a linear transformation carrying
a unit ball to $E(d,l)$, we see that twice the logarithm of the volume of the ellipsoid differs by a constant from
\beq\label{eq:defg}
g(d) := g(d,l) := n \ln f(d,l) - \ln \det ADA^T = n \ln f(d,l) + p(d),
\eeq
where $p(d) := - \ln \det ADA^T$ is a standard barrier function designed to keep $ADA^T$ positive definite.
(Here and below, we define $- \ln \det M$ to be $+ \infty$ if $M$ is not positive definite, even if its determinant is positive.)
We could therefore consider an algorithm that iterates values of $d$ or of $(d,l)$ to minimize $g$, or maybe its upper bound
\beq\label{eq:deftg}
\tilde g(d) := \tilde g(d,l) := n f(d,l) + p(d) - n,
\eeq
which has the advantage of being defined even if $f$ is nonpositive, and is a convex function of $d$.
However, while we have expressions for these functions and their derivatives, they involve the inverse or determinant
of $ADA^T$, which is costly to evaluate for each new $d$. It therefore makes sense to consider coordinate descent
algorithms, changing just a single component of $d$ at each iteration, since then $ADA^T$ is modified by a rank-one update
and its inverse (or Cholesky factorization) and determinant are simple to update. It turns out that such coordinate descent
algorithms are exactly variants of the ellipsoid algorithm, as we shall see in Section \ref{sec:decsteps}.

\section{The ``best'' lower bound}\label{sec:bestlb}

Recall that Section 3 of Burrell and Todd \cite{bt85} shows that a lower bound on a constraint function can be generated by 
any ellipsoid represented as above (with $f(d,l) = 1$). Thus, given an index $j$, we can calculate
\beq\label{eq:deflam}
\lambda := \gamma_j D (A^T (\bar y - \gamma_j^{-1} B a_j) - r),
\eeq
where $\gamma_j := (a_j^T B a_j)^{1/2}$, which satisfies $A \lambda = - a_j$, and then
\[
\theta(\lambda) := l^T \lambda_- - u^T \lambda_+
\]
provides a lower bound on $a_j^T y$ over $P$. This $\lambda$ can be converted into a nonnegative
\[
\hat \lambda := \Lambda \lambda_- + \lambda_+
\]
which also satisfies $A \hat \lambda = - a_j$ and $\theta(\hat \lambda) = - u^T \hat \lambda = \theta(\lambda)$, equations
which directly certify the lower bound from $A^T y \leq u$. We call such $\lambda$'s dual vectors, because they certify 
lower bounds via duality.

In both the SEA and the OEA, a variety of ellipsoids is considered at each iteration. In the SEA, there is the original ellipsoid, the 
ellipsoid obtained by decreasing $d_j$ to zero, and the final ellipsoid. In the OEA, there is the original ellipsoid, the ellipsoid obtained
by decreasing $l_j$ until the center satisfies $a_j^T \bar y = u_j$, and the final ellipsoid. All these ellipsoids differ from the
original merely by a different $d_j$ and a different $l_j$. In this section, we obtain the best lower bound of the form
$\theta(\lambda)$ that can be obtained from a class of $\lambda$'s including all those generated as above.

Note that $\lambda$ above is a linear combination of $D(A^T \by - r)$, $DA^T B a_j$, and $e_j$ (with a zero weight
on the latter). We first show how to get an improved lower bound based on a dual vector of the same form, but with
$j$th component zero, as long as the lower bound from $\lambda$ improves on $l_j$ and $a_j^T \by > r_j$.

Note that Proposition 4.2 in \cite{bt85} shows that $d_j \gamma_j^2 < 1$, and so $\lambda_j = d_j \gamma_j (a_j^T \by - r_j  - \gamma_j)
> -1$. Hence
$$
\tilde \lambda := \frac{1}{1 + \lambda_j} (\lambda - \lambda_j e_j) 
$$
is well defined. It also has $j$th component zero, and satisfies $A \tilde \lambda = - a_j$. If $\lambda_j \geq 0$ then we find
$$
\theta(\tilde \lambda) =  \frac{1}{1 + \lambda_j}  \theta(\lambda) + \frac{\lambda_j}{1 + \lambda_j} u_j,
$$
a convex combination of the bound given by $\lambda$ and $u_j$. If $\theta(\lambda) > u_j$,
indicating infeasibility, then $\theta(\tilde \lambda)$ is also greater than $u_j$, and we can use it to
generate a certificate of infeasibility. If $\theta(\lambda) \leq u_j$, then $\theta(\tilde \lambda)$ provides at least as good a
lower bound as $\theta(\lambda)$.

On the other hand, if $\lambda_j < 0$, then we have
$$
\theta(\tilde \lambda) =  \frac{1}{1 + \lambda_j}  \theta(\lambda) + \frac{\lambda_j}{1 + \lambda_j} l_j,
$$
or
$$
\theta(\lambda) = (1 + \lambda_j) \theta(\tilde \lambda) + (-\lambda_j) l_j,
$$
a convex combination of $\theta(\tilde \lambda)$ and $l_j$. In both cases, we either obtain a certificate of
infeasibility, and terminate, or
$$
\theta(\tilde \lambda) \geq \theta(\lambda),
$$
and note that $\tilde \lambda$ is also a linear
combination of $D(A^T \by - r)$, $DA^T B a_j$, and $e_j$.
 
In the SEA, $a_j^T \by > r_j$, i.e., $t_j$ is positive, in the ellipsoid at the start of the iteration because we choose $j$ as a constraint
where $\by$ violates the upper bound; $t_j$ remains positive in the ellipsoid after dropping $a_j$ 
by (23) in \cite{bt85}; and it is still positive for the final ellipsoid since (29) of \cite{bt85} shows that it is multiplied by 
$1 - \hat \sigma > 0$ from its previous value. In the OEA, $t_j$ is positive in the ellipsoid at the start of
the iteration by the choice of $j$; it remains positive after the decrease of $l_j$ by (47) in \cite{LFT}; and it is still positive
for the final ellipsoid by (41) in \cite{LFT}. Thus in all cases we can move from $\lambda$ to $\tilde \lambda$ and either still have a certificate of infeasibility
or obtain at least as good a lower bound. So we now confine ourselves to dual vectors that are linear combinations of $D(A^T \by - r)$, $DA^T B a_j$, and $e_j$
and have $j$th component zero.

Let $d$ and $l$ denote the vectors used in the ellipsoid at the start of the iteration, and let 
$\hat d$ and $\hat l$ correspond to another of the ellipsoids considered in the previous paragraph,
so that they differ from $d$ and $l$ only in their $j$th components. Let $\hat D$, $\hat r$, and $\hat y$
correspond to these new vectors. Then
 $$
 A \hat DA^T = A D A^T + \alpha a_j a_j^T, \qquad (A \hat DA^T)^{-1} = (A D A^T)^{-1} + \beta (A  D A^T)^{-1} a_j a_j^T (A D A^T)^{-1}
 $$
 for some $\alpha$, $\beta$. Next $A \hat D \hat r = A  D  r + \delta a_j$ for some $\delta$, so that
 $$
 \hat y = (A \hat D A^T)^{-1} A \hat D \hat r = \by + \epsilon (A  D A^T)^{-1} a_j
 $$
 for some $\epsilon$,
 and so $\hat D( A^T \hat y - \hat r)$ is a linear combination of $ D (A^T  y -  r)$, $ D A^T B a_j$, and $e_j$.
 
 Similarly, $(A \hat D A^T)^{-1} a_j = \zeta (A  D A^T)^{-1} a_j$ for some $\zeta$, so that
 $\hat D A^T (A  \hat D A^T)^{-1} a_j$ is a linear combination of $D A^T B a_j$ and $e_j$. It follows that
 $\hat \lambda$, like $ \lambda$, is a linear combination of $ D (A^T \by -  r)$, $D A^T B a_j$, and $e_j$.
 
 We have shown that any lower bound generated by one of the ellipsoids we have been considering arises from a dual vector
 which is a linear combination of these three vectors, and that moreover, we can restrict our attention to those dual vectors
 whose $j$th component vanishes.
 
 Let us therefore consider a generic such dual vector
 $$
 \lambda = \mu  D (A^T \by -  r) + \nu  D A^T B a_j + \pi e_j
 $$
 with $A \lambda = -a_j$ and $\lambda_j = 0$. The first condition yields
 $$
 \nu + \pi = -1,
 $$
 while the second gives 
 $$
 \mu  d_j \bar t_j + \nu d_j  \gamma_j^2 + \pi = 0.
 $$
 
 We can now solve the first equation for $\nu$ in terms of $\pi$, substitute in the second, and then solve
 for $\pi$ in terms of $\mu$, and thus express $\lambda$ as a linear
 function of $\mu$. Next we find and sort the $m$ values of $\mu$ where a component of $\lambda$
 vanishes, and then by moving through the sorted values we find $\mu$ to maximize the piecewise-linear concave function 
 $\theta(\lambda)$. This gives us the desired best lower bound from any of the set of ellipsoids under consideration.
 The work for this last phase is $O(m \ln m)$, while that of the remaining computations (in particular, of calculating
 $A^T B a_j$) is $O(mn)$.
 
 (Note: the final $\lambda$ has zero $j$th component, but this may not be true of $\bar \lambda := \Lambda \lambda_- + \lambda_+$.
 It may then be possible to improve the lower bound further, as in moving from $\lambda$ to $\tilde \lambda$ above; 
 note that $\bar \lambda$ is not of the form above so such an improvement does not contradict our argument.
 This is also why we put ``best'' in quotes above. It would be possible to compute the best lower bound achievable from a dual vector
 as above after the improvement just discussed, but this would require a one-dimensional search for the maximum of a piecewise-rational
 function and would require up to twice as much arithmetical work per iteration, and so we did not pursue it.)

\section{A modified potential function}\label{sec:potfn}

Let us define
\[
z(A,u) := \min_{y \in \R^m} \max_i (a_i^T y - u_i)  , \qquad \tau(A,u) := | z(A,u) |.
\]
We use $\tau(A,u)$ as a condition number for the problem (strictly, its inverse is a condition number). If $P$ is empty, then every point has some $a_i^T y - u_i$ positive, and
since $z(A,u)$ can be written as the optimal value of a bounded linear programming problem this shows that $z(A,u)$ is positive.
If $P$ is nonempty, and has positive volume, then there is a point satisfying all constraints strictly, so that $z(A,u)$ is negative.
If all constraint normals $a_i$ are normalized, then $\tau(A,u)$ is the minimum distance each constraint must be relaxed to make the problem feasible 
in the first case ($P$ empty), and the radius of the largest ball contained in the feasible region in the second ($P$ nonempty). (In \cite{LFT}, all constraint normals
were normalized, but that is not needed in what follows.)

We use $\gamma_i(d,l)$ as the semi-width of the ellipsoid $E(d,l)$ in the direction $a_i$:
\beq\label{eq:defgam}
\gamma_i(d,l) := \max\{ a_i^T (y - \by): y \in E(d,l) \} = \sqrt{ f(d,l) a_i^T B a_i }.
\eeq

If $P$ is empty and we have an ellipsoid $E(d,l)$ where the index $j$ with maximal (positive) $a_i^T \by - u_i$
has $\gamma_j(d,l) < \tau(A,u)$, it is not hard to see (Proposition 5.2 of \cite{LFT}) that we can construct from a dual vector
certifying a lower bound for $a_j^T y$ a certificate of infeasibility.  Thus at every iteration, the OEA of \cite{LFT} ensures a decrease of
not only the volume of $E(d,l)$, but also a measure that forces an aggregate decrease in the $\gamma_i$'s. In fact, \cite{LFT} requires that
the weight vector $d$ is strictly positive, and uses the following upper bounds on the $\gamma_i$'s (Proposition 7.1 of \cite{LFT}):
\[
\gamma_i(d,l) \leq \left( \frac{d_i}{f(d,l)} \right) ^{-1/2}.
\]
Then \cite{LFT} uses the potential function
\[
\phi(d,l) := \prod _{i=1}^m \max \left\{ \left( \frac{d_i}{f(d,l)} \right) ^{-1/2} , \frac{m}{m+1} \tau(A,u) \right\}.
\]
Clearly, $\phi$ is bounded below, by $(m \tau(A,u) / [m+1])^m$, and if we decrease it to this value, then all $\gamma_i(d,l)$'s are
below $\tau(A,u)$ and we can prove infeasibility. Thus, \cite{LFT} proves that $\phi$ decreases by a fixed fraction at every iteration, and
hence establishes convergence in the infeasible case.

One disadvantage of using this function is that it forces all iterates to have $d$ strictly positive. Here, we suggest a modified potential function
that avoids this restriction, and we give geometric interpretations to both potential functions (and incidentally to the infeasibility criterion $f(d,l) < 0$).
We simply use the $\gamma_i$'s directly in the potential function instead of their upper bounds:
\beq\label{eq:defpsi}
\psi(d,l) :=  \prod _{i=1}^m \max \left\{ \gamma_i(d,l) , \frac{m}{m+1} \tau(A,u) \right\}.
\eeq

In order to prove that this modified potential function decreases suitably at each iteration, it suffices to prove the following analog of
Lemma 7.1 of \cite{LFT}:

\BL \label{lm:potdec} 
Let $d  \geq 0$ and $l \in \R^m$ satisfy $f(d,l) > 0$, and similarly let $\tilde{d} \geq 0$ and $\tilde{l} \in \R^m$ satisfy $f(\tilde{d},\tilde{l}) > 0$. Let $1 \leq j \leq m$, and suppose that $d$, $l$, $\tilde d$, $\tilde l$ satisfy:  
$$\frac{1}{f(\tilde{d},\tilde{l})} \tilde{d} = \alpha \left(\frac{1}{f(d,l)} d+ \frac{2}{m-1} \frac{1}{\gamma_j(d,l)^2} e_j \right) \ , $$
for a scalar $\alpha \geq \frac{m^2-1}{m^2}$.  If $\gamma_j(d,l) \geq \tau(A,u)$, then 
$$\psi(\tilde{d}, \tilde{l}) \leq \exp \left(-\frac{1}{2(m+1)}\right) \psi(d,l). $$
\EL

\proof
The proof follows that of Lemma 7.1 of \cite{LFT}  (which proves a similar result for $\phi$ instead of $\psi$), using
\[
\nu_i(d) := \max \left\{ \gamma_i(d,l) , \frac{m}{m+1} \tau(A,u) \right\}
\]
instead of
\[
\mu_i(d) := \max \left\{ \left( \frac{d_i}{f(d,l)} \right) ^{-1/2} , \frac{m}{m+1} \tau(A,u) \right\}.
\]
That proof is quite technical, and its details are not important to the rest of this paper. We therefore just highlight the differences when using the modified potential
function $\psi$. First note that
(with $B$ as usual denoting $(ADA^T)^{-1}$)
\[
\frac{1}{f(\tilde{d},\tilde{l})} A \tilde D A^T = \alpha \frac{1}{f(d,l)} \left( A D A^T + \frac{2}{m-1}\frac{1}{a_j^T  B a_j}  a_j a_j^T \right)
\]
so that
\beq\label{eq:adaupdate} 
f(\tilde{d},\tilde{l}) (A \tilde D A^T)^{-1} = \frac{1}{\alpha} f(d,l) \left( B - \frac{2}{m+1}\frac{1}{a_j^T  B a_j}  Ba_j a_j^T B \right).
\eeq
It follows that
\[
\gamma_j(\tilde d, \tilde l)^2 =  \frac{1}{\alpha} \frac{m-1}{m+1} \gamma_j(d,l)^2 \leq \left(\frac{m}{m+1}\right)^2  \gamma_j(d,l)^2
\]
using the bound on $\alpha$, and similarly
\[
\gamma_i(\tilde d, \tilde l)^2 \leq  \frac{1}{\alpha} \gamma_i(d,l)^2 \leq \frac{m^2}{m^2 -1} \gamma_i(d,l)^2
\]
for $i \neq j$. With a few extra arguments to take care of the maxima in the definition of $\psi$, these inequalities yield
the desired result as in the proof of Lemma 7.1 in \cite{LFT}, using the inequality
\[
\frac{m}{m+1}  \left( \frac{m^2}{m^2 -1} \right) ^{\frac{m-1}{2}} \leq \exp \left(-\frac{1}{2(m+1)}\right).
\]
\qed

We conclude this section by giving some geometric intuition regarding these potential functions, or rather their
versions without taking the maxima,
\[
\hat \phi(d,l) := \prod _{i=1}^m  \left( \frac{d_i}{f(d,l)} \right) ^{-1/2}
\]
and
\beq\label{eq:defpsihat}
\hat \psi(d,l) :=  \prod _{i=1}^m  \gamma_i(d,l) .
\eeq
As the product of $m$ terms, these seem to be related to volumes of $m$-dimensional objects. This is true for the
first, but not quite for the second.

Recall that $E(d,l)$ is an ellipsoid in $\R^n$ containing all feasible $y$'s. Let us view this in the space of the
slacks $s := A^T y \in \R^m$, assuming $d > 0$. Note that a feasible $s$ satisfies $r - v \leq s \leq r + v$, so lies in the ellipsoid
\[
E_s^1 := \{ s \in \R^m: (s - r)^T D (s - r) \leq v^T D v \}.
\]
The volume of this $m$-dimensional ellipsoid is the volume of the $m$-dimensional ball times
\[
\prod_{i=1}^m \left( \frac{d_i}{v^T D v} \right) ^{-1/2},
\]
which differs from $\hat \phi(d,l)$ by replacing $f(d,l) = v^T D v - \bar t ^T D \bar t$ (from (\ref{eq:deff2})) by $v^T D v$.
However, note that we are interested in slacks of the form $A^T y$, so we can take a slice of this ellipsoid of the form
\[
\{ A^T y \in \R^m: (A^T y - (A^T \by - \bar t))^T D (A^T y - (A^T \by - \bar t)) \leq v^T D v \}.
\]
But since $A D \bar t = 0$, this quadratic inequality simplifies to $(A^T y - A^T \by)^T D (A^T y - A^T \by) \leq v^T D v - \bar t^T D \bar t$,
so that all feasible slack vectors $s = A^T y$ lie in the ellipsoid
\[
E_s^2 := \{ s \in \R^m: (s - A^T \by)^T D (s - A^T \by) \leq f(d,l) \},
\]
whose volume is exactly $\hat \phi(d,l)$ times that of the $m$-dimensional ball.

We remark that these two $m$-dimensional ellipsoids, $E_s^1$ centered at $r$ and the similar but smaller one $E_s^2$
centered at $A^T \by$, are analogous to a globe and the smaller globe enclosing all points of a given latitude.
We also see that $f(d,l) < 0$ corresponds to the case that the subspace $\{ A^T y \}$ completely misses the
ellipsoid centered at $r$.
Finally, if $d$ has some zero components, then the two $m$-dimensional ellipsoids become ellipsoidal cylinders, of
infinite volume, while if $A D A^T$ is positive definite, their intersection with the subspace $\{ A^T y \}$ is an
$n$-dimensional ellipsoid embedded in $\R^m$.

The other (perturbed) potential function, $\hat \psi(d,l)$, does not appear to be the volume of an $m$-dimensional object,
but as the product of the semi-widths $\gamma_i(d,l)$, it can be viewed as an $m$-dimensional measure
of the $n$-dimensional ellipsoid $E(d,l)$.

\section{The ellipsoid method as coordinate descent and decrease steps}\label{sec:decsteps}

Let us suppose we have some nonnegative $d$ with $ADA^T$ positive definite and some lower bounds
$l$ certified by $\Lambda$ as in (\ref{eq:Lambda}). We can then define $\by$, and we stop if this is feasible.
So assume not. Then if $f(d,l) \leq 0$, we can construct a certificate of infeasibility as in Section \ref{sec:elliprep}. In fact,
we may as well also check whether $D \bar t$ provides a certificate of infeasibility, and then stop; this may happen even if
$f(d,l)$ is positive.

We therefore assume that $\by$ is not feasible, and that $f(d,l) > 0$, so we can scale $d$ to make $f(d,l)$ equal to 1.
We now have an ellipsoid $E := E(d,l)$ that contains $P$.
Standard ellipsoid methods proceed as follows. First, an index $j$ is chosen so that $\by$ violates the $j$th constraint:
$a_j^T \by > u_j$. Then the new ellipsoid is chosen by adjusting $d_j$ and possibly $l_j$. We have seen how to choose the ``best'' 
possible value for $l_j$. In this section we concentrate on updating $d_j$ in order to decrease $g$ in (\ref{eq:defg}) (and/or possibly $\phi$ or $\psi$).
Since we are changing a single component of $d$, this can be viewed as coordinate descent, and then we may want to consider
{\em decreasing} $d_j$ as well as increasing it if this leads to good progress. It turns out that this is possible if $\by$ satisfies the 
$j$th constraint handily.

\subsection{Decrease steps}

Assume we have chosen a particular index $j$ and we have scaled $d$ so that $f(d,l) = 1$.
For simplicity, we write $\gamma_j$ for $\gamma_j(d,l)$. We consider the implications of updating $d$ to
\beq\label{eq:dplus}
d_+ := d_+(\sigma) := d + \frac{\sigma}{(1-\sigma) \gamma_j^2} e_j,
\eeq
with $\sigma < 1$.
This form is chosen to mesh with previous work: if $0 \leq \sigma < 1$, the resulting ellipsoid is defined by
a quadratic inequality that is a linear combination
of that defining the old ellipsoid (with a weight of $1-\sigma$), and that defining the slice
$l_j \leq a_j^T y \leq u_j$ (with a weight of $\sigma / \gamma_j^2$). Now we consider also the possibility of choosing a negative value of $\sigma$.
Of course, in this case the containment of the intersection of the current ellipsoid and the slice within the
new ellipsoid is not guaranteed, but the containment of the feasible region is guaranteed as long as $d$
remains nonnegative.

If we choose a negative value of $\sigma$, there are three particular values we need to consider:
one that reduces $d_j$ to zero, one that makes $f(d,l)$ equal to zero, and one that minimizes $g(d,l)$.
Note that as $\sigma$ goes from 0 to $-\infty$, $\sigma/[(1 - \sigma) \gamma_j^2]$ goes from 0 to $- \gamma_j^{-2}$, but this is no
great limitation, since by Proposition 4.1 of \cite{bt85}, $d_j \leq \gamma_j^{-2}$, with equality only if making $d_j$ zero makes
$ADA^T$ singular; moreover, the proposition assures that $A D_+ A^T$ is positive definite.

Let us determine the effects on $\by$,
$f$ and $p$ of making such a change, where $\sigma < 1$ can be positive or
negative. This result is similar to those obtained in \cite{tomv} and \cite{bt85}, and
a restatement of Proposition D.1 of \cite{LFT} in our notation.

\BP\label{pr:update}
Suppose $f(d,l) = 1$ and $d_+$ is as given in (\ref{eq:dplus}) above. Let 
\beq\label{eq:alphabeta}
\alpha := \frac{a_j^T \by - u_j}{ \gamma_j} = \frac{\bar t_j - v_j}{ \gamma_j}, \qquad \beta := \frac{a_j^T \by - l_j} {\gamma_j} = \frac{\bar t_j + v_j}{\gamma_j}.
\eeq
Then
\beq\label{eq:sigchanges1}
B_+ := (AD_+A^T)^{-1} = B - \sigma \frac{Ba_j a_j^T B}{a_j^T B a_j}, \quad \by_+ := B_+ AD_+r = \by - \sigma \frac{\alpha + \beta}{2 \gamma_j}  Ba_j, 
\eeq
and
\beq\label{eq:sigchanges2}
\zeta(\sigma) := f(d_+,l) = 1 - \alpha \beta \sigma + \frac{(\beta - \alpha)^2}{4} \frac{\sigma^2}{1 - \sigma} , \quad \pi(\sigma)) := p(d_+) = p(d) + \ln(1 - \sigma).
\eeq
\EP
Note that $\alpha$ and $\beta$ are convenient measures for the depths of the cuts. The $j$th constraint can be written as
$\alpha \leq - a_j^T (y - \by) / \gamma_j \leq \beta$. We always have $\alpha < \beta$, while $\alpha > 0$ signifies
that $\by$ violates the $j$th (upper-bound) constraint. If $\alpha > 1$ then the $j$th constraint fails to intersect the ellipsoid, and a certificate of
infeasibility can be constructed. If $\alpha \leq -1$, then all the points in the ellipsoid satisfy the $j$th constraint. Similar statements apply to $\beta$ with
respect to the $j$th lower-bound constraint.

\vspace{.1in}
\proof
First note that $AD_+A^T = ADA^T + [\sigma/(1-\sigma)] a_j a_j^T / a_j^T B a_j$, so the formula for $B_+$ follows from the rank-one update formula, as
does the equation $ \det B_+ =  (1 - \sigma)  \det B$, which leads to the formula for $\pi(\sigma)$.

Next,
\beann
\by_+ & = & (B - \frac{\sigma}{\gamma_j^2} B a_j a_j^T B) (A D r + \frac{\sigma}{(1 - \sigma) \gamma_j^2} r_j a_j)  \\
    & = & \by - \sigma \gamma_j^{-2} B a_j ( a_j^T \by - \frac{1}{1 - \sigma} r_j + \frac{\sigma}{1 - \sigma} r_j)  \\
    & = & \by - \sigma \gamma_j^{-2} (a_j^T \by - r_j) B a_j = \by - \sigma \gamma_j^{-2} \bar t_j B a_j =\by - \sigma \frac{\alpha + \beta}{2 \gamma_j} B a_j.
\eeann

Then
\beann
\by_+^T (AD_+A^T) \by_+ & = & \by_+^T A D_+ r \\
  & = & (\by - \sigma \gamma_j^{-2} \bar t_j Ba_j)^T (ADr + \frac{\sigma}{(1 - \sigma ) \gamma_j^2} r_j a_j) \\
  & = & \by^T(ADA^T)\by - \sigma\gamma_j^{-2} \bar t_j a_j^T \by + \frac{\sigma}{(1 - \sigma) \gamma_j^2} r_j a_j^T \by
                     - \frac{\sigma^2}{(1 - \sigma) \gamma_j^2} r_j \bar t_j  \\
  & = & \by^T(ADA^T)\by + \frac{\sigma}{(1 - \sigma) \gamma^2}(- [1 - \sigma] \bar t_j (\bar t_j + r_j) + r_j (\bar t_j + r_j) - \sigma r_j \bar t_j) \\
  & = & \by^T(ADA^T)\by + \frac{\sigma}{(1 - \sigma) \gamma_j^2}(r_j^2 - \bar t_j^2 + \sigma \bar t_j^2) =  
           \by^T(ADA^T)\by + \frac{\sigma}{(1 - \sigma) \gamma_j^2} r_j^2 - \sigma \gamma_j^{-2} \bar t_j^2.
\eeann
Finally,
\[
l^T D_+ u = l^T (D + \frac{\sigma}{(1 - \sigma) \gamma_j^2} e_j e_j^T) u = l^TDu +  \frac{\sigma}{(1 - \sigma) \gamma_j^2} l_j u_j =
            l^TDu +  \frac{\sigma}{(1 - \sigma) \gamma_j^2}(r_j^2 - v_j^2),
\]
so that
\beann
\zeta(\sigma) & = & \by_+^T (AD_+ A^T) \by_+ - l^T D_+ u \\
   & = & \by^T (ADA^T) \by - l^T D u +  \frac{\sigma}{1 - \sigma} \gamma_j^{-2} v_j^2 - \sigma \gamma_j^{-2} \bar t_j^2 = 
       1 - \sigma \gamma_j^{-2} (\bar t_j^2 - v_j^2) +  \frac{\sigma^2}{1 - \sigma} \gamma_j^{-2} v_j^2.
\eeann
This gives the desired result on noting that $\gamma_j^{-2} (\bar t_j^2 - v_j^2)  = \gamma_j^{-1} (\bar t_j - v_j) \gamma_j^{-1} (\bar t_j + v_j) = 
\alpha \beta$ and $\gamma_j^{-1} v_j = (\beta - \alpha) / 2$.
\squareforqed

\subsection{Critical values of $\sigma$}

Now we identify the three critical values of $\sigma$. First, the value at which the $j$th component of $d_+(\sigma)$ hits zero is
\[
\sigma_0 := - \frac{d_j}{1 - d_j \gamma_j^2}.
\]

Second, the value (if it exists) where $\zeta(\sigma)$ hits zero, by the proposition above, is a root of the quadratic
\[
[(\beta - \alpha)^2 + 4 \alpha \beta] \sigma^2 - 4 (1 + \alpha \beta) \sigma + 4 = 0,
\]
or
\[
 \frac{2 (1 + \alpha \beta) \pm \sqrt{4 (1+\alpha\beta)^2 - 4 (\alpha + \beta)^2}} {(\alpha + \beta)^2} = 
  2 \frac{1+\alpha\beta \pm \sqrt{(1-\alpha^2)(1-\beta^2)}}{(\alpha+\beta)^2}
\]
if $\alpha + \beta \neq 0$, and
\[
 (1 + \alpha \beta)^{-1} = (1 - \alpha^2)^{-1} = (1 - \beta^2)^{-1}
\]
if $\alpha + \beta = 0$. Note that, if $\alpha + \beta = 0$ and $\alpha \geq -1$, then this last expression is undefined or greater than 1
and is not a valid value for $\sigma$. Similarly, if $\alpha + \beta \neq 0$, we must have $|\alpha|$ and $|\beta|$ both at most 1 or
both greater than one for the square root to exist. In the former case $1 + \alpha \beta > 0$ and by the arithmetic-geometric
mean inequality, the smaller root is at least $2 [1 + \alpha \beta - (1 - \alpha^2/2 - \beta^2/2)] / (\alpha + \beta)^2 = 1$.
Thus the only time that a meaningful value of $\sigma$
can give a zero value of $\zeta$ is when $\alpha < -1$ and $\beta > 1$, and then the root closest to zero is
\[
\sigma_{\zeta} := \left\{ \ba{ccl}  2 \frac{1+\alpha\beta + \sqrt{(1-\alpha^2)(1-\beta^2)}}{(\alpha+\beta)^2} & \mbox{if} & \alpha + \beta \neq 0, \\
    (1 - \alpha^2)^{-1} & \mbox{if} & \alpha + \beta  = 0 \ea \right].
 \]
 
 The last critical value is where $\eta(\sigma) := g(d_+(\sigma))$ attains its minimum. Note that $\eta(\sigma) = n \ln \zeta(\sigma) + \pi(\sigma)$,
 so that 
 \[
 \eta'(\sigma) = \frac{1}{(1 - \sigma)\zeta(\sigma)} \left[ n (1 - \sigma) \left( - \alpha \beta  +
    \frac{(\beta - \alpha)^2}{4}\frac{2\sigma - \sigma^2}{(1 - \sigma)^2} \right)
   - \left(1 - \alpha \beta \sigma + \frac{(\beta - \alpha)^2}{4} \frac{\sigma^2}{1 - \sigma}\right) \right].
 \]
 We see that $\eta'(0)$ is negative if $\alpha \beta > -1/n$ (so that $\sigma$ should be increased to a positive value) and
 positive if $\alpha \beta < -1/n$ (so that $\sigma$ should be decreased to a negative value. We are interested in values of $\sigma$ less than 1
 with $\zeta(\sigma)$ positive,
 and then $\eta'(\sigma)$ is a positive multiple ($1/[4(1-\sigma)^2 \zeta(\sigma)]$) of the quadratic
 \[
 - (n+1)(\alpha + \beta)^2 \sigma^2 + (2n (\alpha + \beta)^2 + 4(1 + \alpha \beta)) \sigma - 4(1 + n \alpha \beta).
 \]
 If $\alpha + \beta = 0$, this is in fact linear, with a root at
 \beq\label{eq:siggam1}
 \sigma_{\eta} := \frac{1 + n \alpha \beta}{1 + \alpha \beta} = \frac{1 - n \beta^2}{1 - \beta^2}.
 \eeq
 Otherwise, if $1 + n \alpha \beta$ is positive, it has two positive roots, and we would like to increase
 $\sigma$ to the smaller, which is
 \beq\label{eq:siggam2}
 \sigma_{\eta} := \frac{2(1+n\alpha\beta)+n(\alpha+\beta)^2-\rho}{(n+1)(\alpha+\beta)^2},
 \eeq
 with
 \[
 \rho := \left( 4(1-\alpha^2)(1-\beta^2)+ n^2(\beta^2 - \alpha^2)^2\right)^{1/2}.
 \]
 These are exactly the formulae given in Todd \cite{tomv}. If $1 + n\alpha\beta$ is negative, the quadratic has one
 positive and one negative root, and we would like to decrease $\sigma$ to the negative root, which is given by exactly
 the same formula $\sigma_{\eta}$ given above.
 
 In the case that $1 + n \alpha \beta$ is positive, the ellipsoid $E(d_+(\sigma_{\eta}))$ is the minimum-volume ellipsoid
 containing a slice of the current ellipsoid $E(d,l)$, A similar statement is true when $1 + n \alpha \beta$ is negative,
 but now involving the two end-pieces of the current ellipsoid.
 
 \BT
 Suppose $-1 \leq \alpha \leq \beta \leq 1$, with $\alpha > -1$ and/or $\beta < 1$, and let
 \[
 \bar E_{\alpha \beta} := \{ y \in E(d,l): - a_j^T y \leq - a_j^T \by + \alpha \gamma_j \, \,
  \mathrm{or} \, - a_j^T y \geq - a_j^T \by + \beta \gamma_j \}.
 \]
 Then if $1 + n \alpha \beta$ is nonnegative, the minimum-volume ellipsoid containing $\bar E_{\alpha \beta}$ is $E(d,l)$, while
 if $1 + n \alpha \beta < 0$, it is $E(d_+(\sigma_{\eta}),l)$.
 \ET
 
 \proof
 We follow exactly the proof technique in \cite{tomv}.  We first transform to a situation where $E(d,l)$ is the unit ball and $- a_j$ becomes
 the first unit vector. Then the very same arguments of  Proposition 1 and Theorem 1 of \cite{tomv} provide lower bounds on the volumes of
 containing ellipsoids, since the $2 n + 1$ or $4n$ points used in those proofs also lie in the (transformed) set $\bar E_{\alpha \beta}$.
 
 Next, we construct ellipsoids achieving these bounds exactly as in Theorem 2 of \cite{tomv}. The only difference is that now, if
 $y$ lies in the (transformed) set $\bar E_{\alpha \beta}$, it satisfies $y^T y \leq 1$ and $(e_1^T y - \alpha)(e_1^T y - \beta) \geq 0$.
 We take $1 - \sigma_{\eta}$ times the first and add $\sigma_{\eta}$ times the second and add. But since the second
 multiplier is negative, it reverses the sense of the second inequality and we get a valid less-than-or-equal-to inequality. The rest of the
 proof is the same.
 \qed

 Our strategy is now clear. If we choose $j$ giving $\alpha \beta > -1/n$, we should increase $\sigma$ to $\sigma_{\eta}$
 as given by (\ref{eq:siggam1}) or (\ref{eq:siggam2}) above. This is like the usual deep-cut ellipsoid method. On the other hand, if
 $\alpha \beta < -1/n$, we should decrease $\sigma$. 	If $\alpha < -1$ and $\beta > 1$, we can compute $\sigma_{\zeta}$.
 If this is at least $\sigma_0$, we should decrease $\sigma$ to $\sigma_{\zeta}$, and the resulting value
 of $f(d,l)$ will be zero, allowing us to deduce the feasibility or infeasibility of the problem. If this does not occur, either because $\alpha \geq 1$
 or $\beta \leq 1$ or $\sigma_{\zeta} < \sigma_0$, we decrease $\sigma$ to the maximum of $\sigma_0$ and $\sigma_{\eta}$
 as given by (\ref{eq:siggam1}) or (\ref{eq:siggam2}) above.
 
 \subsection{Guaranteed reduction of volume but not potential}
 
 We now give conditions under which a decrease step can guarantee a reduction in the volume of the ellipsoid.
 
 \BT\label{th:decsteps}
 Suppose $n \geq 2$. Then, if
 \[
 \alpha \beta \leq - \frac{2}{n} \quad \mbox{and  } \max\{\alpha, -\beta \} \leq - \frac{2}{n},
 \]
 there is some $\hat \sigma < 0$ with
 \[
 \tilde g(d_+(\hat \sigma),l) \leq \tilde g(d,l) - \frac{1}{4n},
  \]
  where $\tilde g$ is defined in (\ref{eq:deftg}).
  Moreover, either
  \itemize
  
  \item[(a)] $d_+(\sigma_{\zeta}) \geq 0$ and $f(d_+(\sigma_{\zeta}),l) \leq 0$; or
  
  \item[(b)] $f(d_+(\sigma_0),l) > 0$ and $\vol(E(d_+(\sigma_0,l) \leq \vol(E(d,l))$; or
  
  \item[(c)] $d_+(\sigma_{\eta}) \geq 0$ with $f(d_+(\sigma_{\eta}),l) > 0$ and 
 \[
 \vol(E(d_+(\sigma_{\eta}),l) \leq \exp\left( - \frac{1}{8n} \right) \vol(E(d,l)).
 \]
 \ET
 
\proof
We will be working with $g(d,l)$, which is twice the logarithm of the volume of $E(d,l)$, and its
upper bound $\tilde g(d,l)$.
Recall that $\tilde g(d,l) = n f(d,l) + p(d) - n$, so that by Proposition \ref{pr:update}
\[
\Delta \tilde g := \tilde g (d_+(\sigma), l) - \tilde g(d,l) = - n \alpha \beta \sigma + n \frac{(\beta - \alpha)^2}{4}\frac{\sigma^2}{1 - \sigma} 
   + \ln (1 - \sigma) \leq - (n \alpha \beta + 1) \sigma + n \frac{(\beta - \alpha)^2}{4} \sigma^2
  \]
 for negative $\sigma$. We now distinguish three cases.
 
 If $\alpha \geq -1$ and $\beta \leq 1$, then $n \alpha \beta + 1 \leq -1$ and the last term on the right-hand side is at most $n \sigma^2$,
 so that the right-hand-side above is at most $\sigma + n \sigma^2$, and by choosing $\hat \sigma = - 1/(2n)$ we see that we
 can achieve the desired decrease.
 
 Next, if $-1 \leq \alpha \leq - 2/n$ and $\beta > 1$ (the argument is analogous if $\alpha < -1$ and $2/n \leq \beta \leq 1$),
 then $n \alpha \beta + 1 \leq - 2 \beta + 1 \leq - \beta$ and $(\beta - \alpha)^2 \leq (\beta + 1)^2 \leq 4 \beta^2$, so
 that the right-hand side is at most $\beta \sigma + n \beta^2 \sigma^2$, and by choosing $\hat \sigma = - 1/(2 n \beta)$
 we see that we can again achieve a decrease of at least $1/(4n)$.
 
 Finally, if $\alpha < -1$ and $\beta > 1$ and $\beta \geq - \alpha$ (the argument is analogous if $-\alpha > \beta$),
 then $n \alpha \beta + 1 \leq - n \beta + 1$ and $(\beta - \alpha)^2 \leq (2 \beta)^2 = 4 \beta^2$, so that the
 right-hand side is at most $(n \beta - 1) \sigma + n \beta^2 \sigma^2$, and by choosing $\hat \sigma = - (n \beta - 1)/(2 n \beta^2)$,
 we see that we can achieve a decrease of at least $(n \beta - 1)^2 / (4 n \beta^2) = (n - 1/\beta)^2 / (4n) \geq (n - 1)^2 / (4n) \geq 1/(4n)$.
 This proves the first part.
 
 Next, if $\sigma_{\zeta}$ is well-defined and at least $\sigma_0$, then (a) holds, and we can establish feasibility or infeasibility.
 So assume not. Then if $\hat \sigma \leq  \sigma_0$, an upper bound on $\tilde g$, and hence one on $g$, is decreasing as
 we move from $\sigma = 0$ down to $\sigma_0$, and it is tight at $\sigma = 0$, so case (b) holds. 
 Lastly, assume $\sigma_{\zeta}$ is not well-defined or is less than $\sigma_0$, and that $\hat \sigma > \sigma_0$.  Then either
 $\sigma_{\eta} \leq \sigma_0$,  so that the volume of $E(d_*(\sigma),l)$ is decreasing as we move from $\sigma = 0$ down to $\sigma_0$,
 and so case (b) holds, or $\sigma_{\eta} > \sigma_0$, and then 
 \[
 g(d_+(\sigma_{\eta}),l) \leq g(d_+(\hat \sigma),l) \leq \tilde g(d_+(\hat \sigma),l) \leq \tilde g(d,l) - \frac{1}{4n} = g(d,l) - \frac{1}{4n},
 \]
 so that case (c) holds. \qed
 
 Unfortunately, we cannot guarantee a suitable decrease in the potential functions $\phi$ and $\psi$ for negative $\sigma$. Indeed, let us suppose all
 terms $\mu_i$ and $\nu_i$ are defined by their first arguments. Then to first order, $f(d_+(\sigma),l)$ is $1 - \alpha \beta \sigma$.
 All components $d_i$, $i \neq j$, are unchanged, while
 \[
 \frac{(d_+(\sigma))_j}{d_j} = 1 + \frac{\sigma}{(1 - \sigma) d_j \gamma_j^2},
 \]
 and since $d_j$ can be arbitrarily close to zero, this ratio can decrease arbitrarily fast as $\sigma$ decreases. This implies that we cannot control $\phi$.
 By Proposition \ref{pr:update}, we know
 \[
 \gamma_j(d_+(\sigma),l)^2 = (1 - \sigma) \gamma_j^2,
 \]
 However, the best we can achieve for $i \neq j$, using the Cauchy-Schwarz inequality, is
 \[
 \gamma_i(d_+(\sigma),l)^2 \leq (1 - \sigma) \gamma_i^2,
 \]
 and since there are $m$ such terms we cannot guarantee a decrease in $\psi$ either. Hence decrease steps do not seem to be
 possible in the OEA while maintaining complexity guarantees for the infeasible case.
 
Suppose we only try decrease steps when the conditions of Theorem \ref{th:decsteps} hold. Then the ellipsoid method 
in the feasible case will enjoy similar theoretical guarantees to those without such steps. Indeed, in a decrease step, 
either case (a) holds and we establish feasibility or infeasibility; or case (b)
holds (we call this a {\em drop} step) and the volume of the ellipsoid does not increase; or case (c) holds, and the volume decreases
by a factor similar to that for a ``normal'' step. Thus we only need to bound the drop steps. But if $d_j$ decreases to zero, it must have been
positive because it was one of the original at most $m$ components of $d$, or because it was increased from $0$ in a previous step (we
call these {\em add} steps). Hence the iteration bound is at most multiplied by eight (four because of the less significant volume reduction,
and two because each drop step must be paired with an earlier increase or add step) and incremented by $m$.
 
 \section{Improving an iteration}\label{sec:improving}
 
 In this section we see how the ideas of the previous sections can be incorporated into the iterations of the SEA and the OEA.
 Note that both of these include some counter-intuitive steps. In the SEA, we first decrease $d_j$ to zero before
 updating the lower bound and then increasing $d_j$ again. In the OEA, we first decrease the lower bound used in the
 ellipsoid before updating the ellipsoid, and the value of $l_j$ used in the ellipsoid is not necessarily the best available.
 Here we remedy these drawbacks and possibly improve the efficiency of the algorithms without sacrificing
 any theoretical guarantees.
 
 Note that most ellipsoid updating steps perform a balancing act: the ellipsoid is squeezed in one direction while
 growing in other directions, so that the parameters must be appropriately chosen to obtain an overall decrease in
 the appropriate measure. This is doubly so for the OEA, since it seeks to simultaneously
 decrease the volume of the ellipsoid and the potential function $\phi$. However, some changes (analogous
 to masking in crowded indoor spaces and getting vaccinated in the time of COVID) are unequivocally beneficial, involving no trade-offs.
 
 \BD A move from $(d,l)$ to $(\tilde d, \tilde l)$ is a {\em Pareto-improving} step if $\tilde d \geq d$ and $f(\tilde d, \tilde l) \leq f(d,l)$,
 with at least one inequality strict.
 \ED
 
 Of course, we are interested in the case that both $l$ and $\tilde l$ are vectors of valid lower bounds, and that $d \geq 0$.
 
 Let us examine the effect of such a step. By (\ref{eq:sigchanges2}), possibly applied repeatedly, $p(d)$ is reduced, so that
 $g(d,l)$ is decreased. By its definition, $\phi(d,l)$ is (maybe not strictly) decreased. And since (\ref{eq:sigchanges1}), possibly applied repeatedly,
 shows that $a_i^T B a_i$ decreases for each $i$, $\psi(d,l)$ is (maybe not strictly) decreased.
 
 Here are some examples of Pareto-improving steps:
 \begin{itemize}
 
 \item[(a)] Increasing $d_j$ if $v_j(d,l) < | \bar t_j(d,l) |$ (i.e., if $a_j^T \by > u_j$ or $a_j^T \by < l_j$), by (42) in \cite{LFT};
 
 \item[(b)] Decreasing $l_j$ if $a_j^T \by > u_j$ by (48) in \cite{LFT};
 
 \item[(c)] Increasing $l_j$ if $a_j^T \by < u_j$, also by (48) in \cite{LFT};
 
 \item[(d)] Increasing both $d_j$ and $l_j$ while keeping $a_j^T \by = u_j$; and
 
 \item[(e)] Taking a minimum-volume ellipsoid updating step if $\delta$ in (3.7) of \cite{tomv} is at most 1, e.g., if
  $\alpha \geq 1/n$ and $\beta = 1$ ($d$ increases while $f(d,l)$ remains 1).
  
  \end{itemize}
  
  The only one of these cases that requires proof is (d), but it is crucial to our improvements. Case (b) corresponds to the first step
  in an iteration of the OEA in \cite{LFT}; case (c), completely analogous, could be applied at the end of such an iteration, either increasing
  $l_j$ to the best certified lower bound or increasing $a_j^T \by$ to $u_j$. Case (e) is included to illustrate that certain minimum-volume
  ellipsoid updates, if the parameters are suitable, involve no trade-offs and will simultaneously decrease the potential functions $\phi$
  and $\psi$. Now let us examine case (d) in detail.

\BP\label{prop:pareto}
Suppose $d \geq 0$, with $ADA^T$ positive definite, and $l$ is a vector of lower bounds for the constraints of $P$. Let
the center $\by$ of $E(d,l)$ satisfy $a_j^T \by = u_j$ where $d_j > 0$, and let $\tilde l_j$, with $l_j < \tilde l_j < u_j$, 
also be a valid lower bound for $a_j^T y$ for $y \in P$.
Define
\[
\mu := \frac{\tilde l_j - l_j}{u_j - \tilde l_j} d_j.
\]
Then, with $\tilde d := d + \mu e_j$, $\tilde l = l + (\tilde l_j - l_j) e_j$, $E(\tilde d, \tilde l)$ also has center $\by$ and the move
from $(d,l)$ to $(\tilde d, \tilde l)$ is a Pareto-improving step.
\EP
(Note that if $d_j = 0$, $\mu = 0$ and $d$ is unchanged. Since the ellipsoid does not depend on the $j$th constraint, it is
unchanged when $\tilde l$ replaces $l$, so that $\by$ and $f(d,l)$ are unchanged and the conclusion is true except that no
inequality in the definition of Pareto-improving steps is strict.)

\vspace{.1in}
\proof
The inequality defining $E(d,l)$ can be written
\beq\label{eq:ineq1}
d_j (a_j^T y - u_j)(a_j^T y - l_j) + \sum_{i \neq j} d_i (a_i^T y - u_i)(a_i^T y - l_i) \leq 0.
\eeq
If we add to this $\mu$ times
\beq\label{eq:ineq2}
(a_j^T y - u_j) (a_j^T y - u_j) \leq 0,
\eeq
we obtain
\[
(d_j + \mu) (a_j^T y - u_j)(a_j^T y - \tilde l_j) +  \sum_{i \neq j} d_i (a_i^T y - u_i)(a_i^T y - l_i) \leq 0,
\]
since
\[
\tilde l_j = \frac{d_j}{d_j + \mu} l_j + \frac{\mu}{d_j + \mu} u_j
\]
from the definition of $\mu$. This is the inequality defining $E(\tilde d, \tilde l)$. Since $\by$ minimizes both the
left-hand side of (\ref{eq:ineq1}) (uniquely) and the left-hand side of (\ref{eq:ineq2}), it also minimizes uniquely the
left-hand side of the final inequality, implying that  it is the center of $E(\tilde d, \tilde l)$.

It remains to show that $f(d,l)$ is not increased. But using the equation (\ref{eq:deff2}), we see that 
\[
f(d,l) = d_j v_j^2 - d_j \bar t_j^2 + \sum_{i \neq j} (d_i v_i^2 - d_i \bar t_i^2),
\]
and  only
the terms indexed by $j$ change. But since $a_j^T \by = u_j$, we find that $v_j = \bar t_j$ both before and after the change,
and thus in both cases, the two terms cancel and hence $f(d,l)$ is unchanged.
\qed

We now outline improvements to the SEA and the OEA. First, we choose the index $j$ as either one that is maximally violated (maximum
$(a_i^T \by - u_i) / \gamma_i(d,l)$ in the SEA or $(a_i^T \by - u_i) / \| a_i \|$ in the OEA), or, for the SEA, an index for which a decrease or drop
step (Section \ref{sec:decsteps}) would be worthwhile. Next, if we choose a violated constraint,
we can use the technique of Section \ref{sec:bestlb} to either generate a lower bound
at least as good as those in previous algorithms, or obtain a certificate of infeasibility. We can then proceed as follows.

For the SEA, instead of decreasing $d_j$ to zero, we first take a Pareto-improving step of type (a) to decrease $a_j^T \by$ to $u_j$. Using (\ref{eq:sigchanges1}),
we see that we should replace $d$ by $d_+(\sigma)$, where
\[
\sigma := 2 \frac{a_j^T \by - u_j}{(\alpha + \beta) \gamma_j}.
\]
Next we can take a Pareto-improving step of type (d) to increase both $d_j$ and $l_j$, increasing the latter to the bound found in the previous paragraph.
Finally, we take a usual ellipsoid-updating step from the current ellipsoid, increasing $d_j$ to minimize the volume of the new ellipsoid.

It is worth noting that exactly the same ellipsoid would have been obtained by first decreasing $d_j$ to zero, then updating $l_j$ to the
``best'' lower bound, and then taking a usual ellipsoid-updating step, since with the same updated $l_j$, there is only one volume-minimizing
$d_j$. However, we find the three-step procedure described above more intuitive, and it paves the way for our improvement of the OEA.

In \cite{LFT}, the first step of the OEA is to decrease $l_j$ until $a_j^T \by = u_j$, and then take an ellipsoid-updating step, finally adjusting $l_j$ to
correspond to the final ellipsoid. The result is that $l_j$ will rarely be the best lower bound found. As an alternative, we can perform the first two steps of the
modified SEA described above, resulting in $l_j$ being the best lower bound found and the equality $a_j^T \by = u_j$. (An alternative, which will lead
to exactly the same ellipsoid, is to decrease $l_j$ until $a_j^T \by = u_j$, and then take a Pareto-improving step of type (d), increasing both
$d_j$ and $l_j$, until $l_j$ reaches the best lower bound, but we prefer the first motivation.) We can then replace $d$ by $d_+(\sigma)$, where
\[
\sigma = \frac{2}{m+1}.
\]
Since $\alpha = 0$ and $\beta \leq 1$ (we use a bound no less than that given by the current ellipsoid)  (\ref{eq:sigchanges1}) gives
\[
f(d_+,l) \leq 1 + \frac{1}{m^2 - 1} = \frac{m^2}{m^2 - 1},
\]
and so Lemma 7.1 of \cite{LFT} and Lemma \ref{lm:potdec} above assure us that we obtain a suitable decrease in $\phi$ and $\psi$.
Indeed, we can do slightly better. Note that $\ln \hat \phi(d_+(\sigma),l)$ is bounded by an expression like 
$g(d_+(\sigma),l)$, but with $m$ replacing $n$. We can therefore use a value
for $\sigma$ that is like $\sigma_{\eta}$ in (\ref{eq:siggam2}), but with $m$ replacing $n$, and then we get a reduction at least as good as that from
$\sigma = 2/(m+1)$ as above. The proof parallels that of Lemma 7.2 in \cite{LFT}.

\section{Computational results}\label{sec:comp}

Here we give the results of some preliminary computational testing.

\subsection{Problem generation}

We consider problems of sixteen different sizes. The number of variables, $n$, is 60, 125, 250, or 500;
the number of inequalities, $m$, is either $1.4$, $2$, $2.8$, or $4$ times $n$ (photographers of a certain age
may recognize these numbers). For each such pair $n$, $m$, we generate ten feasible problems and ten infeasible problems
as follows.

For feasible problems, we first generate $A$ as an $n \times m$ matrix with independent standard Gaussian entries. We then
generate an $n$-vector $y_0$ whose entries are each 100 times a standard Gaussian random variable, again all independent. Finally,
we set $u$ as $A^T y_0$ plus a vector of ones, so that $y_0$ is feasible.

For infeasible problems, we start by generating $A$ and $y_0$ as above. We then generate a vector $x$ of $m$ independent uniform
random variables in $[0,1]$, and replace $A$ by $A - (1 / e^T x) A x e^T$, with $e$ an $n$-vector of ones, so that $A x = 0$.
We set $u$ equal to $A y_0$ plus an $n$-vector of independent standard Gaussian random variables, and then replace $u$ by $-u$ if
$u^T x > 0$. Then $x$ certifies the infeasibility of the system $A^T y \leq u$, but the (now infeasible) region is again somewhat centered about $y_0$ (or
the similarly distributed $- y_0$).

In either case, we apply our algorithms to the system $A^T y \leq u$, without knowing whether or not it is feasible, hoping to generate either a feasible
solution or a certificate of infeasibility.

\subsection{Obtaining the initial system}

We investigated three ways to convert a problem generated as above to a form suitable for our algorithms, hence with
bounds on the variables. 

\subsubsection{Big $M$ initialization}

The simplest method is just to augment the constraints $A ^T y \leq u$ with bounds $ - M e \leq y \leq M e$, where $e$ is an $n$-vector of ones
and $M$ is a large constant. While this method has no theoretical justification, we include it to compare our other methods to, and also because it serves as the basis for the two-phase method below. We typically use $M = 10,000$ which is large compared to the expected size of a component of $y$, of the order
of $100$ for problems generated as above. Hence we do not expect to render feasible problems infeasible by adding these bounds. Similarly,
if we generate a certificate of infeasibility, we expect it to involve only the original constraints of the system, not the added bounds. These
expectations held true in our experiments. We also tried other values of $M$ to determine how it affected the performance of the algorithms. 

\subsubsection{Freund-{Vera} initialization}

The second method is to introduce an additional variable to homogenize the problem, as in Freund and Vera \cite{FandV3} . 
Thus we replace $A^T y \leq u$ by
\[
A^T y - u \eta \leq 0, \qquad - \eta \leq 0.
\]
(We actually want a solution with $\eta$ positive, so that $y / \eta$ will solve the original system.)
With this homogeneous system, we can add arbitrary bounds, so we require $-e \leq y \leq e$, with
$e$ as above, and $\eta \leq 1$ (we already have a lower bound of $0$ on $\eta$).
Note that the system above is always feasible, with $y$ and $\eta$ both zero. If we happen to hit this
solution, we do not terminate, but continue, pretending the constraint $ - \eta \leq 0$ is violated.
If we find a feasible solution with $\eta$ positive, we have a feasible solution to the original system.
We could also find a weak infeasibility certificate, a vector satisfying all the requirements except that the strict inequality is satisfied
weakly, at equality. Then, since all the added upper bounds ($y \leq e, \, -y \leq e, \, \eta \leq 1$) have right-hand sides 1, the weights on all of
these must be zero. It follows that we have $x$, $\xi$ satisfying
\[
A x = 0, \quad - u^T x- \xi = 0, \quad x \geq 0, \, \xi \geq 0.
\]
If in fact $\xi$ is positive, then $A x = 0$, $u^T x = - \xi < 0$, and $x \geq 0$, so we have
a certificate of infeasibility for the original system. (In addition, the system above with $\xi > 0$ is the alternative system for the 
homogeneous system with $ - \eta < 0$.) If not, at least we have a weak infeasibility certificate.

Use of a homogeneous system is very attractive in eliminating the need to add artificial bounds and determine appropriate values
to be used, but the approach also has disadvantages. If the original system is infeasible, then the corresponding homogeneous
system is on the boundary between feasible and infeasible systems, so its condition number is infinite. More practically, it
is clearly necessary computationally to choose tolerances very carefully, to estimate whether an inexact solution is close
enough to feasibility and has $\eta$ sufficiently positive to claim we have found a feasible solution to the original system,
and to estimate whether an inexact infeasibility certificate has $\xi$ sufficiently positive to claim we have determined
infeasibility of the original system. This caused very few problems for the algorithms where decrease of weights $d_j$ was allowed,
but if not, it was obviously hard to generate certificates with zero or close to zero weights on the added bounds, since these
are obtained from Theorem \ref{th:negrhs} or from lower bound weights and usually involve constraints with positive $d_j$'s. Hence a number of 
infeasible instances mistakenly were judged feasible.

\subsubsection{Two-phase method}

Finally, we discuss a technique that is based on the fact that we really want to let our big $M$ tend to infinity, and the observation
that the algorithms are quite insensitive to its size. What happens as $M$ grows to infinity? Scaling and then taking the limit,
we are led to the phase-1 problem
\[
A^T y \leq 0, \qquad -e \leq y \leq e.
\]
We can apply the SEA or the OEA to this system. If we find a point with $A^T y < 0$, we can scale it to find a point satisfying $A^T y \leq u$,
our original system. If instead, $A^T y \leq 0$ but with equality in at least one component, we choose such an index $j$ and continue
the iterations as if the constraint were strictly violated. Finally, we may obtain a weak certificate of infeasibility, that is, a vector $(x,\tilde x, \hat x)$
satisfying $(x, \tilde x, \hat x) \geq 0, \, Ax + \tilde x - \hat x = 0, \, e^T \tilde x + e^T \hat x \leq 0$.  
Then it is easy to see that the weights $(\tilde x, \hat x)$ on the added bounds are zero, so that
$Ax = 0, \, x \geq 0$. 
However, this does not give a weak certificate of infeasibility for the original system, because $u$ has been ignored.
By examining the different cases of infeasibility certificates in Sections \ref{sec:elliprep} and \ref{sec:bestlb}, we see that $x_j$ being zero usually implies that
$d_j$ is zero, unless either an unlikely coincidence occurs, like $a_j^T \bar y = r_j$ in Theorem \ref{th:negrhs}, or $x_j$ is zero because of our
search of the piecewise-linear function $\theta(\lambda)$ in Section \ref{sec:bestlb}. In the latter cases, if $d_j$ is positive and $j$ corresponds to a bound,
we ignore the certificate and continue the iterations. However, if all positive $d_j$'s correspond to original constraints with $x_j$ positive, then we can use this
$d$ to generate a starting ellipsoid for the original system. Indeed, we can use $Ax = 0, \, x \geq 0$ to generate lower bounds $l_j$ for every $j$ with
$x_j$ positive. We then have upper and lower bounds for every constraint with $d_j$ positive, and this gives our ellipsoid $E(d,l)$.
We do not have lower bounds for the remaining constraints, but these can be generated when such a constraint is chosen as violated.
We now move to phase 2, applying the algorithm to the original system starting with $E(d,l)$. Notice that again we have to be careful in 
setting our tolerances in order to recognize weak certificates of infeasibility, and we might occasionally obtain false indications of (in)feasibility due to
numerical inaccuracies.

For phase 1 to terminate, we need to be able to reduce the weights $d_j$ on the added bounds to zero or to negligible values compared
to the rest. This is
very hard for algorithms that never decrease weights, and so we only use this method for variants of the SEA that allow decrease and drop steps.

\subsection{Experimental results}

We ran several versions of our algorithms to investigate the separate effects of our new lower bounds and of allowing decrease and drop 
steps. We give detailed results for three versions of the SEA and one version of the OEA, and make some comments about our
other findings.

For the OEA, to maintain the theoretical guarantees of \cite{LFT}, we need to choose the violated constraint $j$ with the largest
$(a_j^T \bar y - u_j) / \|a_j\|$. Having chosen $j$, we perform the iteration as described at the end of Section \ref{sec:improving}, including the
acceleration in the last paragraph. For the SEA, we compute $j_{\max}$ and $j_{\min}$ maximizing and minimizing $(a_j^T \bar y - u_j) / \gamma_j$, 
respectively, with the second search confined to those  $j$'s with positive $d_j$. 
We first consider $j = j_{\min}$. If the resulting $\alpha_{\min}$ and $\beta_{\min}$
from (\ref{eq:alphabeta}) satisfy $\alpha \beta \leq -2/n$ and dropping $j$ would keep $A D A^T$ positive definite and not increase the 
volume of the ellipsoid, we perform a drop step with this $j$. We prioritize drop steps to remove redundant constraints and, in the case of
the big $M$ or 2-phase methods, to remove bounds. If we cannot use a drop step, we calculate $\alpha_{\max}$ and $\beta_{\max}$ from (\ref{eq:alphabeta})
for $j = j_{\max}$ and then choose an add/increase step for $j_{\max}$ or a decrease step for $j_{\min}$ according as
$\min(1,\alpha_{\max}) \min(1,\beta_{\max})$ is further from or closer to $-1 / n$ than $\max( -1,\alpha_{\min}) \min(1, \beta_{\min} )$.  Having chosen $j$,
we complete the iteration as described in Section \ref{sec:improving}.

The average numbers of iterations for the SEA are given in Table 1. 
Several observations can be made. For all the initialization methods, the lightly constrained feasible problems were the easiest to solve,
followed by the infeasible problems and then the more highly constrained feasible problems. For some reason, for the big-$M$ and two-phase
methods, the feasible problems with $m = 2n$  were harder to solve than the more highly constrained problems. For infeasible problems, all
methods were comparable, with the number of iterations very highly correlated with the dimension $n$ but seemingly independent of the
number $m$ of constraints. For feasible problems, the methods differed quite a bit, but there was no clear winner.
These observations are confirmed by fitting power laws to the data. For example, for the Freund-Vera initialization
on feasible problems, the best fit was $.08 m^{1.72}$, while for infeasible problems it was $.19 n^{1.78}$. The combination of deep cuts, 
improved lower bounds, and decrease and drop steps seems to have dropped the exponent from the worst case quadratic level to between
$1.7$ and $1.8$, but this is nowhere near the linear rate necessary to compete with pivoting (or interior-point) methods. To highlight this, Table 1 
also includes the average number of iterations required on the same problems by the {\tt linprog} routine in MATLAB,
using the dual simplex option and solving $\max \{f^T y : A^T y \leq u \}$ with $f$ identically zero. These are better by a considerable factor, particularly
on the larger problems, and the iterations are no more costly.

When we modified the SEA to use the lower bound suggested in \cite{bt85}, the results varied from 3\% faster to 50\% slower; the geometric
mean of the ratios was 1.20. The results are given in Table 2. 
Next we used the best bound, but eliminated decrease and drop steps. For this comparison, we restricted ourselves
to the Freund-Vera and big $M$ methods. The increase-only algorithm varied from 7\% faster to 249\% slower, with a geometric mean of the ratios
of 1.93. The full results appear in Table 3. 
We conclude that the improved lower bound had only a slightly beneficial effect, but that decrease and drop stops are very advantageous.

We also investigated the effect of the value of $M$ in the big $M$ method. We only looked at the problems with $n$ equal to 60, 125, or 250.
Since theory suggests that the number of iterations might grow linearly with $\log M$, we might expect a growth of 33\% if we increase $M$ from 1000
to 10,000, and 25\% more if we further increase it to $100,000$. However,  the number of iterations for $M = 10,000$ ranged from 2\% faster to
15\% slower compared to that for $M = 1,000$ for feasible problems, with a geometric mean of the ratios of 1.06, while for $M = 100,000$, the algorithm
ranged from 1\% faster to 16\% slower, with the geometric mean of the ratios still 1.06. For the infeasible problems, the effect was even smaller: all
the numbers were within 3\%. We conclude that the effect of $M$ is marginal, at least when using the SEA with its decrease and drop steps.

The average numbers of iterations for the OEA are given in Table 4. 
We see that these are distinctly worse than those for the SEA, although given the compromises made 
in the original OEA we feel that the fact that results are of the same order as for the SEA shows the value of the improvements described in 
Section \ref{sec:improving}. Some runs terminated with numerical problems: either the lower-triangular Cholesky factor of $ADA^T$ became close to singular,
or $ADA^T$ itself was judged not sufficiently positive definite during a refactorization, called for when inaccuracies in updated quantities became
noticeable. Here we make two comparisons. The first is with the increase-only variant of the SEA of the previous paragraph.
We find that the OEA results are between 7\% faster and 207\% slower (discounting the runs that terminated with numerical problems), with
the geometric mean of the ratios being 1.97. Secondly, we wondered whether part of the poor performance of the OEA was due to its selecting
the index $j$ based on scaling by the norm of $a_j$, rather than by its ellipsoidal norm $\gamma_j$, which would lead to a deeper cut. So we changed
to selecting $j$ in the latter way. Compared to this variant, the OEA results varied from 5\% faster to 18\% slower (ignoring the runs that terminated 
with numerical problems or that ran into the iteration limit). We conclude tentatively (due to the numerical problems) that choosing $j$ to
maintain theoretical guarantees has only a small effect.

Overall, we feel we have given ellipsoid algorithms for linear inequalities that iterate a sequence of containing ellipsoids, choosing one index
$j$ and only updating the corresponding weight $d_j$ and lower bound $l_j$, the very best chance to show their potential, and they have come up wanting.
Even the most sophisticated coordinate-descent method must suffer compared to a gradient- or higher-order-based algorithm, and the speed of the
resulting iterations due to simpler linear algebra is unable to compensate. In high dimensions, the decrease in the volume of the ellipsoid
(and of the potential function) is just too slow to be competitive. 
One might hope that the sequence of ellipsoids generated would resemble a balloon blown up and
then released, rushing round the room and rapidly reducing its volume; instead, one is left with the image of a large soap bubble blown to amuse a child,
which vibrates charmingly in one direction and then another, but seems not to get smaller -- until it pops!

\begin{table}[h!] \label{comp2}
\begin{center}
\begin{tabular}{ | r | r | r | r | r | r | r | r |  r | r | r | r |}
\hline
\multicolumn{2}{| c |}{SEA} & \multicolumn{2}{| c |}{Freund-Vera} & \multicolumn{2}{| c |}{Big $M$} & \multicolumn{2}{| c |}{Two-Phase}  & \multicolumn{2}{| c |}{Dual Simplex}  \\
 \multicolumn{2}{| c |} { } & \multicolumn{2}{ |c |}{Initialization} & \multicolumn{2}{| c |}{Initialization} & \multicolumn{2}{| c |}{Method} & \multicolumn{2}{| c |}{Algorithm}  \\
\cline{1-10} 
 $n$ & $m$ & Feasible & Infeasible & Feasible & Infeasible &Feasible & Infeasible & Feasible & Infeasible \\
\hline
\hline
 60 & 84 &  168.1 & 294.4 &  223.4 & 293.4 & 230.1 & 298.6  & 44.2 &  73.4  \\
\hline
60 & 120 & 448.7 &    283.0       &     589.2       &     283.5        &   587.0       &     283.5     & 62.0 &  83.9 \\
\hline
60 & 168 & 575.1  & 291.7  & 569.7  &  290.1 & 422.5  &  289.9   &  73.3 & 83.2   \\
\hline
60 & 240 &  574.6  & 298.4  & 587.3  &   302.3    &  426.3  &   301.3  &  83.2 & 86.3   \\
\hline
125 & 175 &  477.5  &  1012.5 & 566.7   &  1029.6  & 565.9  & 1044.9  &  103.2 &  165.2      \\
\hline 
125 & 250 & 1690.2 & 1020.2 & 2076.9 & 1017.2 & 2142.4 & 1011.7  &  146.4 & 190.3   \\
\hline
125 & 350 & 2334.0 & 1031.3 & 1648.3 & 1039.3 & 1384.2 & 1032.8  & 172.5  & 196.2  \\
\hline
125 & 500 & 2209.8 & 1082.3 & 1661.7 & 1079.4 & 1415.1 & 1080.4  & 206.1 & 206.4  \\
\hline
250 & 350 & 1143.8 & 3571.5 &  1224.3 &  3551.1 & 1256.1 & 3567.1 &  236.4 & 360.9   \\
\hline
250 & 500 & 7216.6 & 3473.3 & 6848.7 & 3468.2 & 7424.4 & 3468.2 & 329.6 & 417.8   \\
\hline
250 & 700 & 8766.8 & 3547.2 & 4819.0 & 3537.8 & 4518.0 & 3545.2 & 397.9 & 446.9   \\
\hline
250 & 1000 & 7647.5 & 3736.0 & 4888.8 & 3728.3 & 4561.9 & 3732.1 & 476.3 & 455.2   \\
\hline
500 & 700 & 2521.9 & 12951.1 & 2648.3 & 12941.7 & 2602.3 & 12977.2 & 535.7 & 769.6   \\
\hline
500 & 1000 & 26069.4 & 12487.0 & 30432.2 & 12460.8 & 31015.1 & 12495.3 & 736.7 &  906.6  \\
\hline
500 & 1400 & 31222.7 & 12627.1 & 15896.6 & 12633.9 & 15805.4 & 12632.8 & 899.8 &  923.8  \\
\hline
500 & 2000 & 30540.2 & 13128.0 & 15958.0 & 13144.6 & 15675.3 & 13139.8 & 1139.3 & 985.4   \\
\hline
\end{tabular}
\end{center}
\caption{Average Number of Iterations for the Standard Ellipsoid and Dual Simplex Algorithms.   } 
\end{table}

\begin{table}[h!] \label{compo}
\begin{center}
\begin{tabular}{ | r | r | r | r | r | r | r | r |  r | r |}
\hline
\multicolumn{2}{| c |}{SEA} & \multicolumn{2}{| c |}{Freund-Vera} & \multicolumn{2}{| c |}{Big $M$} & \multicolumn{2}{| c |}{Two-Phase}  \\
 \multicolumn{2}{| c |} {old lower bound } & \multicolumn{2}{ |c |}{Initialization} & \multicolumn{2}{| c |}{Initialization} & \multicolumn{2}{| c |}{Method}  \\
\cline{1-8} 
 $n$ & $m$ & Feasible & Infeasible & Feasible & Infeasible &Feasible & Infeasible  \\
\hline
\hline
 60 & 84 &  171.6 & 360.0 &  233.6 & 353.1 & 229.7 & 355.5  \\
\hline
60 & 120 & 463.4 &    339.9       &     629.6       &     346.4        &   625.2       &     338.7    \\
\hline
60 & 168 & 575.3  & 339.7 & 599.9  &  357.6 & 463.1  &  343.0     \\
\hline
60 & 240 &  569.0  &  344.4 & 663.6 & 377.3      &  492.6  &   344.7   \\
\hline
125 & 175 &  483.0  &  1359.4 & 590.3 &  1385.5  & 578.5  &   1362.0     \\
\hline 
125 & 250 & 1820.0 & 1291.5 & 2275.5 & 1322.1 & 2395.7  & 1298.7   \\
\hline
125 & 350 & 2380.6 & 1278.2 & 1698.8 & 1279.6 & 1626.5 & 1278.9   \\
\hline
125 & 500 & 2259.7 & 1279.2 & 1860.1 & 1320.7 & 1637.6  & 1289.5 \\
\hline
250 & 350 & 1148.5 & 5166.8 &  1232.1 &  5172.2 & 1246.6 & 5148.0    \\
\hline
250 & 500 & 7528.9 & 4852.7 & 7820.1 & 4852.9 & 8010.1 & 4859.5   \\
\hline
250 & 700 & 9128.6 & 4700.3 & 5615.5 & 4933.0 & 5684.3 & 4708.4   \\
\hline
250 & 1000 & 7937.7 & 4676.2 & 5780.1 & 4774.2 & 5601.8 & 4684.6  \\
\hline
500 & 700 & 2569.5 & 20209.1 & 2657.5 & 20247.1 & 2640.4 & 20189.3   \\
\hline
500 & 1000 & 27541.3 & 18714.4 & 34961.1 & 18701.1 & 35913.6 & 18716.2  \\
\hline
500 & 1400 & 32721.9 & 18029.6 & 20234.0 & 18094.9 & 20845.0 & 18014.4  \\
\hline
500 & 2000 & 29761.3 & 17714.2 & 20008.0 & 17853.5 & 19897.8 & 17704.0   \\
\hline
\end{tabular}
\end{center}
\caption{Average Number of Iterations for the Standard Ellipsoid Algorithm Using the Old Lower Bound.   } 
\end{table}

\begin{table}[h!] \label{compi}
\begin{center}
\begin{tabular}{ | r | r | r | r | r | r |  }
\hline
\multicolumn{2}{| c |}{SEA} & \multicolumn{2}{| c |}{Freund-Vera} & \multicolumn{2}{| c |}{Big $M$}   \\
 \multicolumn{2}{| c |} {  no decrease steps }  & \multicolumn{2}{| c |}{Initialization} & \multicolumn{2}{| c |}{Initialization}  \\
\cline{1-6} 
 $n$ & $m$ & Feasible & Infeasible & Feasible & Infeasible \\
\hline
\hline
 60 & 84 &  156.9 & 955.7 &  220.8 & 756.3   \\
\hline
60 & 120 & 426.6 &    960.1       &     1260.0       &     726.2     \\
\hline
60 & 168 & 615.5  & 1018.9  & 1102.7  &  746.4    \\
\hline
60 & 240 &  628.6  & 1031.7  & 1083.0&   774.9      \\
\hline
125 & 175 &  452.6  &  3138.4 & 538.8  &  2610.8       \\
\hline 
125 & 250 & 1668.5 &  3128.2* & 5073.5 & 2450.8  \\
\hline
125 & 350 & 2513.4 & 3087.0 & 3700.0 & 2446.9  \\
\hline
125 & 500 & 3095.0 & 3246.3 & 3513.6 & 2522.7 \\
\hline
250 & 350 & 1077.3 & 9943.2 &  1152.6 &  8767.8    \\
\hline
250 & 500 & 6988.6 & 9636.5 & 18739.8 & 7982.8    \\
\hline
250 & 700 & 10801.3 & 9843.1 & 12225.5 & 7888.9 \\
\hline
250 & 1000 & 12477.0 & 9771.8 & 11154.2 & 8063.0  \\
\hline
500 & 700 & 2366.2 & 33079.6* & 2499.9 & 30014.4 \\
\hline
500 & 1000 & 25796.9 & 30651.4* & 79708.0 & 26702.4 \\
\hline
500 & 1400 & 51327.2 & 30728.0 & 40222.4 & 26011.5  \\
\hline
500 & 2000 & 49914.1 & 30015.3 & 36197.5 & 26281.5  \\
\hline
\end{tabular}
\end{center}
\caption{Average Number of Iterations for the Standard Ellipsoid Algorithm with No Decrease Steps.  
Starred entries gave one or two false indications of feasibility. } 
\end{table}

\begin{table}[h!]\label{comp3}
\begin{center}
\begin{tabular}{ | r | r | r | r | r | r | }
\hline
\multicolumn{2}{| c |}{OEA} & \multicolumn{2}{| c |}{Freund-Vera} & \multicolumn{2}{| c |}{Big $M$}   \\
 \multicolumn{2}{| c |} { } & \multicolumn{2}{ |c |}{Initialization} & \multicolumn{2}{| c |}{Initialization} \\
\cline{1-6}
 $n$ & $m$ & Feasible & Infeasible & Feasible & Infeasible  \\
\hline
\hline
 60 & 84 &  155.3 & 1938.9 & 302.8 &  1530.6  \\
\hline
60 & 120 & 545.6 &  1920.7$^\dagger$ &  2417.4       &     1513.3             \\
\hline
60 & 168 & 621.6  & ***** & 2223.5  & 1661.5    \\
\hline
60 & 240 &  623.8  & 2424.6 & 2331.7  & 1936.9     \\
\hline
125 & 175 &  476.9  &  6927.4 & 681.3 & 5752.7      \\
\hline 
125 & 250 & 2640.8 & 6775.7 & 9939.3 & 5553.1 \\
\hline
125 & 350 & 3212.9 & 7190.6 & 8119.6 & 5943.0  \\
\hline
125 & 500 & 3632.7 & 8170.7 & 8334.9 & 6900.5  \\
\hline
250 & 350 & 1060.6 &  24583.9$^\dagger$ & 1335.9 &  20656.7  \\
\hline
250 & 500 & 10656.3 & 22527.6 & 35910.7 & 19230.9  \\
\hline
250 & 700 & 13255.5 & 23553.7 & 28108.1 & 20433.9  \\
\hline
250 & 1000 & 14967.7 & 26418.2 & 28531.3 & 23466.8  \\
\hline
500 & 700 & 2220.8 & ***** & 2625.5 & 74527.0  \\
\hline
500 & 1000 & 39625.3 & 76003.3$^\dagger$ & $ 89640.3 $  & 67401.3  \\
\hline
500 & 1400 & 53417.8 & 79077.7$^\dagger$ & 98419.0 & 70930.6  \\
\hline
500 & 2000 & 61351.9 & 88321.6 & 98719.9 & 80721.8 \\
\hline
\end{tabular}
\end{center}
\caption{Average Number of Iterations for the Oblivious Ellipsoid Algorithm. 
$*$ At least one run terminated for numerical reasons. 
$\dagger$ One or two runs gave false indications of feasibility. } 
\end{table}

\bigskip


 {\bf Acknowledgement}
 The author would like to thank  Jourdain Lamperski and Rob Freund 
 for very helpful conversations, and  two anonymous referees for their detailed and constructive comments.
 
 \pagebreak

\bibliographystyle{amsplain}
\bibliography{ellipredux65}

\providecommand{\bysame}{\leavevmode\hbox to3em{\hrulefill}\thinspace}
\providecommand{\MR}{\relax\ifhmode\unskip\space\fi MR }
\providecommand{\MRhref}[2]{%
  \href{http://www.ams.org/mathscinet-getitem?mr=#1}{#2}
}
\providecommand{\href}[2]{#2}
\begin{thebibliography}{10}

\bibitem{bgt}
R.~Bland, D.~Goldfarb, and M.~J. Todd, \emph{The ellipsoid method: a survey},
  Operations Research \textbf{29} (1981), no.~6, 1039--1091.

\bibitem{blo}
J.V. Burke, A.S. Lewis, and M.L. Overton, \emph{The speed of {S}hor's
  r-algorithm}, IMA Journal of Numerical Analysis \textbf{28} (2008), 711?720.

\bibitem{bt85}
B.~P. Burrell and M.~J. Todd, \emph{The ellipsoid method generates dual
  variables}, Mathematics of Operations Research \textbf{10} (1985), no.~4,
  527--715.

\bibitem{ek}
J.G. Ecker and M.~Kupferschmid, \emph{A computational comparison of the
  ellipsoid algorithm with several nonlinear programming algorithms}, SIAM
  Journal on Control and Optimization \textbf{23} (1985), no.~5, 657--674.

\bibitem{FandV3}
Robert~M. Freund and Jorge~R. Vera, \emph{Condition-based complexity of convex
  optimization in conic linear form via the ellipsoid algorithm}, SIAM Journal
  on Optimization \textbf{10} (1999), no.~1, 155--176.

\bibitem{glspaper}
M.~Gr\"{o}tschel, L.~Lov\'{a}sz, and A.~Schrijver, \emph{The ellipsoid method
  and its consequences in combinatorial optimization}, Combinatorica \textbf{1}
  (1981), 169--197.

\bibitem{gls}
\bysame, \emph{Geometric algorithms and combinatorial optimization}, second
  ed., Springer-Verlag, Berlin, 1994.

\bibitem{kp}
R.M. Karp and C.H. Papadimitriou, \emph{On linear characterizations of
  combinatorial optimization problems}, SIAM Journal on Computing \textbf{11}
  (1982), 620--632.

\bibitem{kha}
L.~G. Khachiyan, \emph{A polynomial algorithm in linear programming}, Soviet
  Math. Dokl. \textbf{20} (1979), no.~1, 191--194.

\bibitem{LFT}
J.~Lamperski, R.M. Freund, and M.J. Todd, \emph{An oblivious ellipsoid
  algorithm for solving a system of (in)feasible linear inequalities},
  Mathematics of Operations Research (2020), to appear, arXiv:1910.03114.

\bibitem{lev}
A.~Iu. Levin, \emph{On an algorithm for the minimization of convex functions},
  Mathematics Doklady \textbf{6} (1965), 286--290.

\bibitem{nemyud77}
A.S. Nemirovski and D.B. Yudin, \emph{Optimization methods adapted to the
  ``significant'' dimension of the problem}, Automation and Remote Control
  \textbf{38} (1977), 513--524.

\bibitem{new}
D.J. Newman, \emph{Location of the maximum on unimodal surfaces}, Journal of
  the Association of Computing Machinery \textbf{12} (1965), 395--398.

\bibitem{pr}
M.W. Padberg and M.R. Rao, \emph{The {R}ussian method for linear programming
  {III}: Bounded integer programming}, Tech. report, New York University,
  Graduate School of Business Administration, 1981, Research Report 81-39.

\bibitem{rn}
A.~Rodomanov and Y.~Nesterov, \emph{Subgradient ellipsoid method for nonsmooth
  convex problems}, Mathematical Programming (2022), to appear.

\bibitem{shor}
N.Z. Shor, \emph{Cut-off method with space extension in convex programming
  problems}, Cybernetics \textbf{13} (1977), 94--96.

\bibitem{Shorg}
N.Z. Shor and V.I. Gershovich, \emph{Family of algorithms for solving convex
  programming problems}, Cybernetics \textbf{15} (1979), 502--507.

\bibitem{tke88}
S.~P. Tarasov, L.G. Khachiyan, and I.~I. Erlikh, \emph{The method of inscribed
  ellipsoids}, Soviet Math. Dokl. \textbf{37} (1988), 226--230.

\bibitem{tikh}
V.~M. Tikhomirov, \emph{The evolution of methods of convex optimization},
  American Mathematical Monthly \textbf{103} (1996), 65--71.

\bibitem{tomv}
M.J. Todd, \emph{Minimum volume ellipsoids containing part of a given
  ellipsoid}, Mathematics of Operations Research \textbf{7} (1980), 253--261.

\bibitem{yn}
D.B. Yudin and A.S. Nemirovski, \emph{Informational complexity and effective
  methods of solution for convex extremal problems}, Matekon: Translations of
  Russian and East European Mathematical Economics \textbf{13} (1977), 25--45.

\end{thebibliography}

\end{document}